\documentclass[11pt]{article}
\usepackage{amsmath,amssymb}
\usepackage{setspace}

\usepackage{graphicx,color}
\definecolor{DarkBlue}{rgb}{0.1,0.1,0.5}
\definecolor{Red}{rgb}{0.9,0.0,0.1}

\def\tr{\text{tr}\,}

\begin{document}


\begin{center}
{\Large{\bf Consistency of restricted maximum likelihood
estimators of principal components}}

\vskip.1in {{\bf Running Title:} Consistency of REML estimators}

 \vskip.2in
Debashis Paul\footnote{Paul's research was partially supported by
the National Science Foundation under Agreement No. DMS-0112069,
as part of a research fellowship from SAMSI.


} ~\textit{and} Jie Peng

\vskip.1in \textit{Department of Statistics, University of
California, Davis}

\date{}
\end{center}

\begin{abstract}

In this paper we consider two closely related problems :
estimation of eigenvalues and eigenfunctions of the covariance
kernel of functional data based on (possibly) irregular
measurements, and the problem of estimating the eigenvalues and
eigenvectors of the covariance matrix for high-dimensional
Gaussian vectors. In Peng and Paul (2007), a restricted maximum
likelihood (REML) approach has been developed to deal with the
first problem. In this paper, we establish consistency and derive
rate of convergence of the REML estimator for the functional data
case, under appropriate smoothness conditions. Moreover, we prove
that when the number of measurements per sample curve is bounded,
under squared-error loss, the rate of convergence of the REML
estimators of eigenfunctions is near-optimal. In the case of
Gaussian vectors, asymptotic consistency and an efficient score
representation of the estimators are obtained under the assumption
that the effective dimension grows at a rate slower than the
sample size. These results are derived through an explicit
utilization of the intrinsic geometry of the parameter space,
which is non-Euclidean. Moreover, the results derived in this
paper suggest an asymptotic equivalence between the inference on
functional data with dense measurements and that of the high
dimensional Gaussian vectors.

\end{abstract}

\section{Introduction}

Analysis of functional data, where the measurements per subject, or
replicate, are taken on a finite interval, has been one of the
growing branches of statistics in recent times. In fields such as
longitudinal data analysis, chemometrics, econometrics, the
functional data analysis viewpoint has been successfully used to
summarize data and gain better understanding of the problems at
hand. The monographs of Ramsay and Silverman (2005) and Ferraty and
Vieu (2006) give detailed accounts of the applications of functional
data approach to various problems in these fields. Depending on how
the individual curves are measured, one can think of two different
scenarios - (i) when the curves are measured on a dense grid; and
(ii) when the measurements are observed on an irregular, and
typically sparse set of points on an interval. The first situation
usually arises when the data are recorded by some automated
instrument, e.g. in chemometrics, where the curves represent the
spectra of certain chemical substances. The second scenario is more
typical in longitudinal studies where the individual curves could
represent the level of concentration of some substance, and the
measurements on the subjects may be taken only at irregular time
points. In the first scenario, i.e., data on a regular grid, as long
as the individual curves are smooth, the measurement noise level is
low, and the grid is dense enough, one can essentially treat the
data to be on a continuum, and employ techniques similar to the ones
used in classical multivariate analysis. For example, Hall and
Hosseini-Nasab (2006) derive stochastic expansions of sample PCA
when the sample curves are noise-free and measured on a continuum.
However, in the second scenario, the irregular nature of the data,
and the presence of measurement noise pose challenges and require a
different treatment. Under such a scenario, data corresponding to
individual subjects can be viewed as partially observed, and
noise-corrupted, independent realizations of an underlying
stochastic process. The estimation of the eigenvalues and
eigenfunctions of a smooth covariance kernel, from sparse, irregular
measurements, has been studied by various authors including James,
Hastie and Sugar (2000), Yao, M\"{u}ller and Wang (2005), and Peng
and Paul (2007), among others.

In Peng and Paul (2007), a \textit{restricted maximum likelihood
(REML)} approach is taken to obtain the estimators. REML estimators
are widely used and studied in statistics. For example, the
usefulness of REML and profile REML estimation has been recently
demonstrated in the context of functional linear mixed effects model
by Antoniadis and Sapatinas (2007).  In Peng and Paul (2007), it is
assumed that the covariance kernel can be well-approximated by a
positive-semidefinite kernel of finite rank $r$ whose eigenfunctions
can be represented by $M (\geq r)$ known orthonormal basis
functions. Thus the basis coefficient matrix $B$ of the approximant
belongs to the \textit{Stiefel manifold} of $M \times r$ matrices
with orthonormal columns. The working assumption of Gaussianity
allows the authors to derive the log-likelihood of the observed data
given the measurement times. Then a Newton-Raphson procedure, that
respects the geometry of the parameter space, is employed to obtain
the estimates by maximizing the log-likelihood. This procedure is
based on the formulation of a general Newton-Raphson scheme on
Stiefel manifold developed in Edelman, Arias and Smith (1998). Peng
and Paul (2007) also derive a computationally efficient approximate
cross-validation score for selecting $M$ and $r$. Through extensive
simulation studies, it is demonstrated that the REML estimator is
much more efficient than an alternative procedure (Yao \textit{et
al.}, 2005) based on local linear smoothing of empirical
covariances. The latter estimator does not naturally reside in the
parameter space, even though it has been proved to achieve the
optimal non-parametric convergence rate in the minimax sense under
$l^2$ loss, under the optimal choice of the bandwidth and when the
number of measurements per curve is bounded (Hall, M\"{u}ller and
Wang, 2006). Also, in most situations, our method outperforms the EM
approach of James \textit{et al.} (2000). Although the latter
estimator also aims to maximize the log-likelihood, it does not
naturally reside in the parameter space either, and thus it does not
utilize its geometry efficiently.

The superior numerical performance of the REML estimator motivates
us to conduct a detailed study of its asymptotic properties. In this
paper, we establish consistency and derive the rate of convergence
(under $l^2$ loss) of the REML estimator when the eigenfunctions
have a certain degree of smoothness, and when a stable and smooth
basis, e.g., the cubic B-spline basis with a pre-determined set of
knots, is used for approximating them. The techniques used to prove
consistency differ from the standard asymptotic analysis tools when
the parameter space is Euclidean.
Specifically, we restrict our attention to small ellipsoids around
zero in the tangent space to establish a mathematically manageable
neighborhood around an ``optimal parameter'' (a good approximation
of the ``true parameter'' within the \textit{model space}). We
derive asymptotic results when the number of measurements per curve
grows sufficiently \textit{slowly} with the sample size (referred as
the \textit{sparse case}).  We also show that for a special scenario
of the \textit{sparse case}, when there is a bounded number of
measurements per curve, the risk of the REML estimator (measured in
squared-error loss) of the eigenfunctions has asymptotically
near-optimal rate (i.e., within a factor of $\log n$ of the optimal
rate) under an appropriate choice of the number of basis functions.

Besides the \textit{sparse case}, we consider two other closely
related problems: (i) the estimation of the eigenvalues and
eigenfunctions of a smooth covariance kernel, from dense, possibly
irregular, measurements (referred as the \textit{dense case}); and
(ii) the estimation of the eigenvalues and eigenvectors of a
high-dimensional covariance matrix (referred as the \textit{matrix
case}). In the \textit{matrix case}, we assume that there is
preliminary information so that the data can be efficiently
approximated in a lower dimensional \textit{known} linear space
whose effective dimension grows at a rate slower than the sample
size $n$. The proofs of the results in all three cases utilize the
intrinsic geometry of the parameter space through a decomposition of
the Kullback-Leibler divergence. However, the \textit{matrix case}
and the \textit{dense case} are more closely related, and the
techniques for proving the results in these cases are different in
certain aspects from the treatment of the \textit{sparse case}, as
described in Sections \ref{sec:functional} and \ref{sec:matrix}.

Moreover, in the \textit{matrix case}, we also derive a
\textit{semiparametric efficient score representation} of the REML
estimator (Theorem 4), that is given in terms of the
\textit{intrinsic Fisher information operator} (note that the
residual term is not necessarily $o_P(n^{-1/2})$). This result is
new, and explicitly quantifies the role of the intrinsic geometry of
the parameter space on the asymptotic behavior of the estimators.
Subsequently, it points to an \textit{asymptotic optimality} of the
REML estimators. Here, asymptotic optimality means achieving the
asymptotic minimax risk under $l^2$ loss within a suitable class of
covariance matrices (kernels). We want to point out that, in the
\textit{matrix case}, the REML estimators coincide with the usual
PCA estimates, i.e., the eigenvalues and eigenvectors of the sample
covariance matrix (Muirhead, 1982). In Paul and Johnstone (2007), a
first order approximation of the PCA estimators is obtained by
matrix perturbation analysis. Our current results show that the
efficient score representation coincides with this approximation,
and thereby gives a geometric interpretation to this. The
asymptotically optimal rate of the $l^2$-risk of the REML estimator
in the \textit{matrix case} follows from this representation and the
lower bound on the minimax rate obtained in Paul and  Johnstone
(2007). Asymptotic properties of high-dimensional PCA under a
similar context have also been studied by Fan, Fan and Lv (2007).
Recently several approaches have been proposed for estimating large
dimensional covariance matrices and their eigenvalues and
eigenvectors under suitable sparsity assumptions on the population
covariance, e.g. Bickel and Levina (2007, 2008) and El Karoui
(2008).

At this point, we would like to highlight the main contributions of
this paper. First, we have established the consistency and derived
the rate of convergence of REML estimators for functional principal
components in two different regimes - the \textit{sparse case} and
the \textit{dense case}. In Hall \textit{et al.} (2006), it is shown
that an estimator of functional principal component based on a local
polynomial approach achieves the optimal nonparametric rate when the
number of measurements per curve is bounded. However, to the best of
our knowledge, no results exist regarding the consistency, or rate
of convergence, of the REML estimators in the functional data
context. Secondly, we have derived an efficient score representation
for sample principal components of high-dimensional, i.i.d. Gaussian
vectors. This involves calculation of the intrinsic Fisher
information operator and its inverse, and along the line we also
provide an independent verification that the REML estimates under a
rank-restricted covariance model are indeed the PCA estimates.
Thirdly, we expect that the current framework can be refined to
establish efficient score representation of the REML estimators of
the functional principal components, and therefore the results
obtained in this paper serve as first steps towards studying the
asymptotic optimality of these estimators. Moreover, results
obtained in this paper suggest an \textit{asymptotic equivalence}
between the inference on functional data with dense measurements and
that of the high dimensional Gaussian vectors. Finally, our work
provides useful techniques for dealing with the analysis of
estimation procedures based on minimization of a loss function (e.g.
MLE, or more generally M-estimators) over a non-Euclidean parameter
space for semiparametric problems. There has been some work on
analysis of maximum likelihood estimators for parametric problems
when the parameter space is non-Euclidean (see e.g. Oller and
Corcuera, 1995). However, there has been very limited work for
non/semi-parametric problems with non-Euclidean parameter space.
Recently, Chen and Bickel (2006) establish semiparametric efficiency
of estimators in ICA (Independent Component Analysis) problems using
a sieve maximum likelihood approach.



The rest of the paper is organized as follows. In Section
\ref{sec:functional}, we present the data model for the functional
principal components, and state the consistency results of the
REML estimators. In Section \ref{sec:matrix}, we describe the
model for high-dimensional Gaussian vectors and derive asymptotic
consistency and an efficient score representation of the
corresponding REML estimators. Section \ref{sec:proof_thm1} is
devoted to giving an overview of the proof of the consistency
result for the functional data case (Theorems 1 and 2). Section
\ref{sec:proof_thm3} gives an outline of the proof of consistency
in the \textit{matrix case} (Theorem 3), in particular emphasizing
the major differences with the proof of Theorem 1. Section
\ref{sec:proof_thm4} is concerned with the proof of the score
representation in the \textit{matrix case} (Theorem 4). Section
\ref{sec:discussion} has a summary of the results and a discussion
on some future works. Technical details are given in the
appendices.

\section{Functional data}\label{sec:functional}

In this section, we start with a description of the functional
principal components analysis, and then make a distinction between
the \textit{sparse case} and the \textit{dense case}.  We then
present the asymptotic results and relevant conditions for
consistency under these two settings.

\subsection{Model}\label{sec:setup}

Suppose that we observe data $Y_i = (Y_{ij})_{j=1}^{m_i}$, at the
design points $T_i = (T_{ij})_{j=1}^{m_i}$, $i=1,\ldots,n$, with
\begin{equation}\label{eq:functional_basic}
Y_{ij} = X_i(T_{ij}) + \sigma \varepsilon_{ij},
\end{equation}
where $\{\varepsilon_{ij}\}$ are i.i.d. $N(0,1)$, $X_i(\cdot)$ are
i.i.d. Gaussian processes on the interval $[0,1]$ (or, more
generally, $[a,b]$ for some $a < b$) with mean 0 and covariance
kernel $\overline{\Sigma}_0(u,v) = \mathbb{E}[X_i(u)X_i(v)]$.
$\overline{\Sigma}_0$ has the spectral decomposition
$$
\overline{\Sigma}_0(u,v) = \sum_{k=1}^\infty \overline{\lambda}_k
\overline{\psi}_k(u) \overline{\psi}_k(v)
$$
where $\{\overline{\psi}_k\}_{k=1}^\infty$ are orthonormal
eigenfunctions and $\overline{\lambda}_1 > \cdots >
\overline{\lambda}_r > \overline{\lambda}_{r+1} \geq \cdots \geq 0$
are the eigenvalues. The assumption that the stochastic process has
mean zero is simply to focus only on the asymptotics of the
estimates of eigenvalues and eigenfunctions of the covariance kernel
(i.e., the functional principal components).

Throughout this paper we assume Gaussianity of the observations. We
want to emphasize that, Gaussianity is more of a working assumption
in deriving the REML estimators. But it plays a less significant
role in asymptotic analysis. For the functional data case, only
place where Gaussianity is used is in the proof of Proposition 3,
and even this can be relaxed by assuming appropriate tail behavior
of the observations. Gaussianity is more crucial in the analysis for
the matrix case. The proofs of Proposition 6 and Theorem 4 depend on
an exponential inequality on the extreme eigenvalues of a Wishart
matrix (based on a result of Davidson and Szarek (2001)), even
though we expect the non-asymptotic bound to hold more generally.

In this paper, we are primarily interested in the situation where
the design points are i.i.d. from a distribution with density $g$
(random design). We shall consider two scenarios, to be referred as
the \textit{sparse case} and the \textit{dense case}, respectively.
The \textit{sparse case} refers to the situation when the number of
measurements, $m_i$, are comparatively small (see {\bf B1}). The
\textit{dense case} refers to the situation where the $m_i$'s are
large so that the design matrix (i.e., the matrix of basis functions
evaluated at the time points) has a concentration property (see {\bf
B1'} and {\bf D}). In the latter case, we also allow for the
possibility that the design is non-random.

Next, we describe the \textit{model space}, to be denoted by ${\cal
M}_{M,r} := {\cal M}_{M,r} (\phi)$, (for $1 \leq r \leq M$) for the
REML estimation procedure. The model space ${\cal M}_{M,r}$ consists
of the class of covariance kernels $C(\cdot,\cdot)$, which have rank
$r$, and whose eigenfunctions are represented in a known orthonormal
basis $\{\phi_k\}_{k=1}^M$ of smooth functions. Furthermore, the
nonzero eigenvalues are all distinct. For example, in Peng and Paul
(2007), $\{\phi_k\}_{k=1}^M$ is taken to be an orthonormalized cubic
$B$-spline basis with equally spaced knots. Thus, the model space
consists of the elements $C(\cdot,\cdot) = \sum_{k=1}^r \lambda_k
\psi_k(\cdot)\psi_k(\cdot)$, where $\lambda_1
> \cdots > \lambda_r
> 0$, and $(\psi_1(\cdot),\ldots,\psi_r(\cdot)) =
(\boldsymbol{\phi}(\cdot))^T B$, where $B$ is an $M \times r$ matrix
satisfying $B^T B = I_r$, and $\boldsymbol{\phi}(\cdot) =
(\phi_1(\cdot),\ldots,\phi_M(\cdot))^T$. Note that we do not assume
that $\overline{\Sigma}_0$ belongs to the \textit{model space}. For
the asymptotic analysis, we only assume that it can be
well-approximated by a member of the \textit{model space} (see
condition {\bf C} and \textit{Lemma 1}). We define the \textit{best
approximation error} of the model as $\inf_{\widetilde C \in {\cal
M}_{M,r}(\phi)}
\parallel \overline{\Sigma}_0 - \widetilde C \parallel_F$, where
$\parallel \cdot \parallel_F$ denotes the Hilbert-Schmidt norm.  A
rank $r$ approximation to $\overline{\Sigma}_0$ in ${\cal M}_{M,r}
(\phi)$ can be defined as
$$
\Sigma_{*0}(u,v) = \sum_{k=1}^r \lambda_{*k} \psi_{*k}(u)
\psi_{*k}(v),
$$
with $\lambda_{*1} > \cdots > \lambda_{*r} > 0$, and
\begin{equation*}
(\psi_{*1}(t),\cdots,\psi_{*r}(t)) = (\boldsymbol{\phi}(t))^T B_*,
\end{equation*}
where $B_*$ is an $M \times r$ matrix satisfying $B_*^T B_* = I_r$.
We refer to $\{(\psi_{*k},\lambda_{*k})\}_{k=1}^r$, or equivalently,
the pair $(B_*,\Lambda_*)$, as an \textit{optimal parameter}, if the
corresponding $\Sigma_{*0}$ is a close approximation to
$\overline{\Sigma}_0$ in the sense that, the approximation error
$\parallel \overline{\Sigma}_0 - \Sigma_{*0}\parallel_F$ has the
same rate (as a function of $M$) as the best approximation error.
Henceforth,  $(B_*,\Lambda_*)$ is used to denote an optimal
parameter.

Observe that, under model (\ref{eq:functional_basic}), $Y_i$ are
independent, and conditionally on $T_i$ they are distributed as
$N_{m_i}(0,\overline{\Sigma}_i)$. Here, the $m_i \times m_i$ matrix
$\overline{\Sigma}_i$ is of the form $\overline{\Sigma}_i =
((\overline{\Sigma}_0(T_{ij},T_{ij'})))_{j,j'=1}^{m_i} + \sigma^2
I_{m_i}$. Then the matrix $\Sigma_{*i} = \Phi_i^T B_* \Lambda_*
B_*^T \Phi_i + \sigma^2 I_{m_i}$ is an approximation to
$\overline{\Sigma}_i$, where $\Phi_i :=
[\boldsymbol{\phi}(T_{i1}):\cdots:\boldsymbol{\phi}(T_{im_i})]$ is
an $M \times m_i$ matrix. We shall use $\Lambda$ to denote
interchangeably the $r \times r$ diagonal matrix
diag$(\lambda_1,\ldots,\lambda_r)$ and the $r\times 1$ vector
$(\lambda_1,\ldots,\lambda_r)^T$. Note that, the parameter
$(B,\Lambda)$ belongs to the parameter space $\Omega := {\cal
S}_{M,r} \otimes \mathbb{R}_+^r$, where ${\cal S}_{M,r} = \{A \in
\mathbb{R}^{M \times r} : A^T A = I_r\}$ is the \textit{Stiefel
manifold} of $M \times r$ matrices with orthonormal columns. For
fixed $r$ and $M$, the \textit{REML estimator} of
$\{(\overline{\psi}_k,\overline{\lambda}_k)\}_{k=1}^r$ is defined as
a minimizer over $\Omega$ of the negative log-likelihood (up to an
additive constant and the scale factor $n$):
\begin{equation}\label{eq:restr_loglike}
L_n(B,\Lambda) = \frac{1}{2n} \sum_{i=1}^n \tr(\Sigma_i^{-1}Y_i
Y_i^T) + \frac{1}{2n} \sum_{i=1}^n \log |\Sigma_i|,
\end{equation}
where $\Sigma_i = \Phi_i^T B \Lambda B^T \Phi_i +  \sigma^2
I_{m_i}$.

\subsection{Consistency}\label{subsec:functional_consistency}

We shall present results on consistency of the REML estimators of
functional principal components in the two different regimes
considered above, namely, the \textit{sparse case} (i.e., when the
number of measurements per curve is ``small'') and the \textit{dense
case} (i.e., when the number of measurements per curve is
``large''). Throughout this paper, we assume that $\sigma^2$ is
known, even though Peng and Paul (2007) provide estimate of
$\sigma^2$ as well. This assumption is primarily to simplify the
exposition. It can be verified that all the consistency results
derived in this paper hold even when $\sigma^2$ is estimated. We
make the following assumptions about the covariance kernel
$\overline{\Sigma}_0$.
\begin{itemize}
\item
[{\bf A1}] The $r$ largest eigenvalues of $\overline{\Sigma}_0$
satisfy, (i) $c_1 \geq \overline{\lambda}_1
> \cdots > \overline{\lambda}_r > \overline{\lambda}_{r+1}$
for some $c_1 < \infty$; (ii) $\max_{1\leq j \leq r}
(\overline{\lambda}_j- \overline{\lambda}_{j+1})^{-1} \leq c_2
<\infty$.
\item
[{\bf A2}] The eigenfunctions $\{\overline{\psi}_k\}_{k=1}^r$ are
four times continuously differentiable and satisfy
$$
\max_{1\leq k \leq r}
\parallel \overline{\psi}_k^{(4)} \parallel_\infty
\leq C_0 ~~~\mbox{for some}~~~ 0 < C_0< \infty.
$$
\end{itemize}

\vskip.15in\noindent{\bf SPARSE case.} In this case, we only
consider the situation when $\sigma^2$ is fixed (i.e., it does not
vary with $n$). We shall first deal with the case when $m_i$'s are
bounded. Then we extend our results to the situation when $m_i$'s
increase slowly with sample size, and are of the same order of
magnitude for all $i$ (condition {\bf B1}). We also assume a
boundedness condition for the random design (condition {\bf B2}).
\begin{itemize}
\item
[{\bf B1}] The number of measurements $m_i$ satisfy $\underline{m}
\leq m_i \leq \overline{m}$ with $4 \leq \underline{m}$ and
$\overline{m}/\underline{m}$ is bounded by some constant $d_2 > 0$.
Also, $\overline{m} = O(n^\kappa)$ for some $\kappa \geq 0$.
\item
[{\bf B2}]
For each $i$, $\{T_{ij}:j=1,\ldots, m_i \}$ are i.i.d. from a
distribution with density $g$, where $g$ satisfies
\begin{equation}\label{eq:g_cond}
c_{g,0} \leq g(x) \leq c_{g,1} ~~\mbox{for all}~~x \in
[0,1],~\mbox{where}~~0 < c_{g,0} \leq c_{g,1} < \infty.
\end{equation}
\end{itemize}
Finally, we have a condition on the $l^2$ error for approximating
the covariance kernel in the model space ${\cal M}_{M,r}$. Define
the \textit{maximal approximation error} for an optimal parameter
$(B_*,\Lambda_*)$ as:
\begin{equation}\label{eq:beta_n_def}
\overline{\beta}_n := \max_{1\leq i \leq n} \frac{1}{m_i} \parallel
\overline{\Sigma}_i - \Sigma_{*i} \parallel_F.
\end{equation}
\begin{itemize}
\item [{\bf C}]
$\overline{m}\overline{\beta}_n = O(\sqrt{\frac{M \log n}{n}})$.
\end{itemize}
If we use orthonormalized cubic $B$-spline basis for representing
the eigenfunctions, then {\bf C} follows from {\bf A1}-{\bf A2} and
{\bf B1}-{\bf B2}, if the covariance kernel is indeed of rank $r$:

\vskip.1in\noindent{\bf Lemma 1 :} \textit{If {\bf A1}-{\bf A2} and
{\bf B1}-{\bf B2} hold,  $\overline{\Sigma}_0$  is of rank $r$, and
we use the orthonormalized cubic $B$-spline basis with equally
spaced knots to represent the eigenfunctions, then {\bf C} holds, if
$M^{-1} (n\overline{m}^2/\log n)^{1/9} = O(1)$.}

\vskip.1in\noindent Proof of Lemma 1 follows from the fact that for
a cubic $B$-spline basis, for sufficiently large $M$, we can choose
$(B_*,\Lambda_*)$ such that (i) $\max_{1\leq k \leq r}
\parallel \overline{\psi}_k - \psi_{*k}
\parallel_\infty = O(M^{-4})$ (by {\bf A1} and {\bf A2}), and
(ii) $\overline{\beta}_n = O(M^{-4})$ (see Appendix A). This implies
that $\parallel \overline{\Sigma}_{0} - \Sigma_{*0}\parallel_F =
O(M^{-4})$. The assumption that the covariance kernel is of finite
rank can be relaxed somewhat by considering the true parameter as a
sequence of covariance kernels $\overline{\Sigma}_{0,n}$ such that
the $(r+1)$-th largest eigenvalue $\overline{\lambda}_{r+1,n}$
decays to zero sufficiently fast. Note that in Lemma 1, the use of
B-spline basis is not essential. The result holds under the choice
of any \textit{stable basis} (i.e., the Gram matrix has a bounded
condition number) with sufficient smoothness.

We now state the main result in the following theorem.

\vskip.1in\noindent{\bf Theorem 1 (sparse case):} \textit{Suppose
that {\bf A1}-{\bf A2}, {\bf B1}-{\bf B2} and {\bf C} hold, and
$\overline{m}$ is bounded. Suppose further that $M$ satisfies
\begin{equation}\label{eq:M_condition_sparse}
M \to \infty, ~~~\mbox{such that}~~~ M^{-1} (n/\log n)^{1/9} = O(1)
~~~\mbox{and}~~~ M = o(\sqrt{n/\log n}), ~~~\mbox{as}~~n\to \infty.
\end{equation}
Then, given $\eta
> 0$, there exists $c_{0,\eta} > 0$ such that for $\alpha_n = c_{0,\eta} \sigma
\sqrt{\frac{\overline{m}^2 M\log n}{n}}$, with probability at least
$1-O(n^{-\eta})$, there is a minimizer $(\widehat B,\widehat
\Lambda)$ of (\ref{eq:restr_loglike}) satisfying
\begin{eqnarray*}
\parallel \widehat B - B_* \parallel_F &\leq&
\alpha_n,\\
\parallel \widehat \Lambda - \Lambda_* \parallel_F &\leq& \alpha_n.
\end{eqnarray*}
Moreover, the corresponding estimate of the covariance kernel, viz.,
$\widehat \Sigma_0(u,v) = \sum_{k=1}^r \widehat \lambda_k \widehat
\psi_k(u) \widehat \psi_k(v)$, satisfies, with probability at least
$1-O(n^{-\eta})$,}
\begin{equation*}
\parallel \widehat \Sigma_0 - \overline{\Sigma}_0 \parallel_F =
O(\alpha_n).
\end{equation*}

\vskip.1in\noindent{\bf Corollary 1:} \textit{Suppose that the
conditions of Theorem 1 hold. Then the best rate of convergence
holds if $M \asymp (n/\log n)^{1/9}$, and the corresponding rate is
given by $\alpha_n \asymp (\log n/n)^{4/9}$. For estimating the
eigenfunctions, this is within a factor of ~$\log n$ of the optimal
rate. The optimal rate over a class ${\cal C}$ of covariance kernels
of rank $r$ satisfying conditions {\bf A1}-{\bf A2}, and the random
design points satisfying conditions {\bf B1}-{\bf B2} (with
$\overline{m}$ bounded), is $n^{-4/9}$.}


\vskip.1in\noindent Notice that, the rate obtained here for the
estimated eigenvalues is not optimal. We expect a parametric rate of
convergence for the latter, which can be achieved by establishing an
efficient score representation of the estimators along the line of
Theorem 4. The following result generalizes Theorem 1 by allowing
for $m_i$'s to slowly increase with $n$, and its proof is
encapsulated in the proof of Theorem 1.

\vskip.1in\noindent{\bf Corollary 2:} \textit{Suppose that, {\bf
A1}-{\bf A2}, {\bf B1}-{\bf B2} and {\bf C} hold. Suppose further
that, $\overline{m}$ and $M$ satisfy
\begin{eqnarray}\label{eq:m_M_condition}
&&(i)~ \overline{m}^4 M \log n = o(n), ~~~~(ii)~\max
\{\overline{m}^3 M^{5/2} (\log n)^2, \overline{m}^{7/2} M^2 (\log
n)^{3/2}\} = o(n),
\nonumber\\
&& ~~(iii)~M^{-1}(n\overline{m}^2/\log n)^{1/9} = O(1), ~~~~(iv)~
\overline{m}^2 M^2 \log n = o(n).
\end{eqnarray}
Then the conclusion of Theorem 1 holds. Also, the best rate is
obtained when $M \asymp (n\overline{m}^2/\log n)^{1/9}$, and the
corresponding $\alpha_n \asymp \overline{m}^{10/9} (\log
n/n)^{4/9}$. }

\vskip.1in\noindent Condition (i) is required to ensure that
$\overline{m}^2 \alpha_n^2 = o(1)$; condition (ii) is needed to
ensure that the upper bound in (\ref{eq:Sigma_diff_upper2}) in Lemma
4 is $o(1)$; condition (iii) ensures that {\bf C} holds; and
finally, condition (iv) is used in proving Lemmas 4, 5 and 6 in
Appendix B. A sufficient condition for (\ref{eq:m_M_condition}) to
hold is that $\overline{m} = O(n^{1/5})$ and $M \asymp
(n\overline{m}^2/\log n)^{1/9}$. Notice that the best rate obtained
in Corollary 2 is not optimal in general. It is near-optimal (up to
a factor of $\log n$ of the optimal rate) only when $\overline{m}$
is bounded above (Theorem 1).

\vskip.15in\noindent{\bf DENSE case.} This case refers to the
scenario where the number of time points per curve is large, such
that $\min_{1\leq i \leq n} m_i \to \infty$ sufficiently fast (see
condition {\bf D} and the corresponding discussion). For
simplicity, we assume further that the number of design points is
the same for all the sample curves, which is not essential for the
validity of the results. Denote this common value by $m$. In terms
of the asymptotic analysis, there is an important distinction
between the \textit{sparse case} and \textit{dense case}. For the
purpose of further exposition and the proof of the result on
consistency of REML estimator in the \textit{dense case}, it is
more convenient to work with the transformed data $\widetilde Y_i
= \Phi_i Y_i$. Let $\Gamma_i = \frac{1}{m} \Phi_i \Sigma_i
\Phi_i^T$ and $R_i = \frac{1}{m} \Phi_i \Phi_i^T$. Then $\Gamma_i
= m R_i B \Lambda B^T R_i + \sigma^2 R_i$. Then, a way of
estimating $\{(\overline{\lambda}_k,\overline{\psi}_k)\}_{k=1}^r$
is by minimizing the negative log-likelihood of the transformed
data:
\begin{equation}\label{eq:restr_loglike_trans_dense}
\widetilde L_n(B,\Lambda) = \frac{1}{2n} \sum_{i=1}^n
\tr(\Gamma_i^{-1}\frac{1}{m}\widetilde Y_i \widetilde Y_i^T) +
\frac{1}{2n} \sum_{i=1}^n \log |\Gamma_i|.
\end{equation}
Notice that, if $R_i$'s are non-singular for all $i$, then by
direct computation, we have that the negative log-likelihoods for
the raw data: (\ref{eq:restr_loglike}) and that of the transformed
data: (\ref{eq:restr_loglike_trans_dense}) differ only by a
constant independent of the parameters $B$ and $\Lambda$. Hence,
on the set $\{R_i~\mbox{are non-singular for all}~i\}$, the
estimators obtained by minimizing (\ref{eq:restr_loglike}) and
(\ref{eq:restr_loglike_trans_dense}) are the same. Assumptions
{\bf B1} and {\bf B2} are now replaced by:
\begin{itemize}
\item
[{\bf B1'}] $m = O(n^{\kappa})$ for some $\kappa > 0$.
\item
[{\bf D}] Given $\eta > 0$, there exist constants
$c_{1,\eta},c_{2,\eta}
> 0$ such that
\begin{equation}\label{eq:R_i_bound}
 \mathbb{P}\left(\max_{1\leq i \leq n} \parallel R_i - I_M
\parallel \leq c_{1,\eta} \sqrt{\frac{\sigma^2}{m  \log n}},~~
\max_{1\leq i \leq n} \parallel B_*^T R_i B_* - I_r \parallel \leq
c_{2,\eta} \frac{\sigma^2}{m \log n}\right) \geq 1- O(n^{-\eta}).
\end{equation}
\end{itemize}
Denote the event described in (\ref{eq:R_i_bound}) by $A_{1,\eta}$.
Note that $A_{1,\eta}$ is defined in terms of $\mathbf{T} :=
\{T_{ij}:j=1,\ldots,m; i=1,\ldots,n\}$ alone. We assume throughout
that $\sigma^2 \leq m$ (note that $\sigma^2/m$ can be viewed as the
signal-to-noise ratio. Therefore, for $n$ large enough, on
$A_{1,\eta}$, $R_i$ is invertible for all $i$. The condition {\bf D}
gives concentration of individual $R_i$'s around the identity matrix
and is discussed in more detail at the end of this section. Finally,
we make an assumption about the maximal approximation error
$\overline{\beta}_n$ defined through (\ref{eq:beta_n_def}) which
differs slightly from the condition {\bf C} in the sparse case.
\begin{itemize}
\item [{\bf C'}] Given $\eta > 0$, there is a constant $c_\eta > 0$
such that $\overline{\beta}_n  \leq c_\eta \frac{\sigma^2}{m}
\sqrt{\frac{M\log n}{n}}$ with probability at least
$1-O(n^{-\eta})$.
\end{itemize}
A result similar to Lemma 1 can be proved to ensure condition {\bf
C'} when a stable basis is used.


\vskip.1in\noindent{\bf Theorem 2 (dense case):} \textit{Suppose
that {\bf A1}-{\bf A2}, {\bf B1'}, {\bf C'} and {\bf D} hold, and $m
\geq \sigma^2 > 0$. Then, given $\eta > 0$, there exists $c_{0,\eta}
> 0$ such that for $\alpha_n = c_{0,\eta} \sigma \sqrt{\frac{M\log n}{n m}}$,
with probability at least $1-O(n^{-\eta})$, there is a minimizer
$(\widehat B,\widehat \Lambda)$ of
(\ref{eq:restr_loglike_trans_dense}) satisfying
\begin{eqnarray*}
\parallel (I_M - B_* B_*^T) (\widehat B - B_*) \parallel_F &\leq& \alpha_n,\\
\parallel B_*^T (\widehat B - B_*) \parallel_F &\leq&
\sqrt{\frac{m}{\sigma^2}} \alpha_n,\\
\parallel \widehat \Lambda - \Lambda_* \parallel_F &\leq&
\sqrt{\frac{m}{\sigma^2}}\alpha_n.
\end{eqnarray*}
Further, the corresponding estimated covariance kernel $\widehat
\Sigma_0(u,v) = \sum_{k=1}^r \widehat \lambda_k \widehat
\psi_k(u)\widehat \psi_k(v)$ satisfies, with probability at least
$1-O(n^{-\eta})$,}
\begin{equation*}
\parallel \widehat \Sigma_0 - \overline{\Sigma}_0 \parallel_F =
O\left(\sqrt{\frac{M\log n}{n}}\right).
\end{equation*}

\vskip.1in\noindent The proof of Theorem 2 requires a slight
refinement of the techniques used in proving Theorem 3 stated in
Section \ref{subsec:matrix_consistency}, making heavy use of
condition {\bf D}. To save space, we omit the proof. Note that the
best rate in Theorem 2 implicitly depends on conditions {\bf C'}
and {\bf D} in a complicated way, which is not the optimal rate
for estimating the principal components. The optimal rate for
$l_2$ risk of the eigenfunctions in this context is conjectured to
be of the order $\max\{(\sigma^2/nm)^{8/9},(1/n)\}$ with the
second term within brackets appearing only when $r > 1$. This can
be verified for the case $r=1$ with a refinement of the proof of
Corollary 1.

\vskip.15in\noindent{\bf Discussion on condition {\bf D}:} We shall
only consider the setting of an uniform design - either fixed, or
random. The condition (\ref{eq:R_i_bound}) clearly requires $m$ to
be sufficiently large, since it gives concentration of individual
$R_i$'s around the identity matrix. To fulfil {\bf D}, we also need
some conditions on the basis functions used. Specifically, we
concentrate on the following classes of basis functions. We assume
that the basis functions are at least 3 times continuously
differentiable.
\begin{itemize}
\item[{\bf E1}] (\textit{Sinusoidal basis}) $\max_{1\leq k \leq M}
\parallel \phi_k \parallel_\infty = O(1)$.

\item[{\bf E2}] (\textit{Spline-type basis}) (i) For any $k \in
\{1,\ldots,M\}$, at most for a bounded number of basis functions
$\phi_l$, supp$(\phi_k) \cap$ supp$(\phi_l)$ is nonempty; (ii)
$\max_{1\leq k \leq M} \parallel \phi_k \parallel_\infty =
O(\sqrt{M})$.
\end{itemize}
One of the key observations in the case of functional data is that,
the eigenfunctions $\{\psi_{*k}\}_{k=1}^r$ of the kernel
$\Sigma_{0*}$ (belonging to the model space) have the same degree of
smoothness as the basis $\{\phi_k\}_{k=1}^M$, and the functions
$\{\psi_{*k}\}_{k=1}^r$ and their derivatives are bounded. Also,
notice that, $B_*^T R_i B_* = ((\frac{1}{m} \sum_{j=1}^m
\psi_{*k}(T_{ij})\psi_{*l}(T_{ij})))_{k,l=1}^r$. Based on these
observations, we present some sufficient conditions for
(\ref{eq:R_i_bound}) to hold under the uniform design and bases of
type {\bf E1} or {\bf E2}. We omit the proof, which uses Bernstein's
inequality (in the random design case) and the \textit{Trapezoidal
rule} (in the fixed design case).

\vskip.1in\noindent{\bf Proposition 1:} \textit{Suppose that the
basis is of type {\bf E1} or {\bf E2}. In the case of random,
uniform design, (\ref{eq:R_i_bound}) is satisfied if $(M \log
n)^2/\sigma^{2} = O(1)$, and $\sqrt{m} \log n /\sigma^{2}= O(1)$. In
the case of fixed, uniform design, (\ref{eq:R_i_bound}) holds (with
probability 1) if $\frac{M^2\log
n}{m}(1+\frac{M^{5/2}}{m})^2/\sigma^{2} = O(1)$, and $\log
n/\sigma^{2}= O(1)$. Moreover, in this setting, if the
eigenfunctions $\{\overline{\psi}_k\}_{k=1}^r$ vanish at the
boundaries, and if the basis functions are chosen so that they also
vanish at the boundaries, it is sufficient that $\frac{M^7\log
n}{m^{7/2}}/\sigma^{2} = O(1)$ and $\frac{\log n}{m}/ \sigma^{2} =
O(1)$.}


\vskip.1in\noindent Note that, two obvious implications of
Proposition 1 are that (i) $m$ needs to be rather large; and (ii)
$\sigma^2$ may need to grow with $n$, in order that {\bf D} holds.

\vskip.1in\noindent {\bf Remark 1:}
It is to be noted that even though the consistency results for the
functional data problem are proved under a specific choice of the
basis for representing the eigenfunctions, viz., the
(orthonormalized) cubic B-spline basis with equally spaced knots,
this is by no means essential. The main features of this basis are
given in terms of the various properties described in Appendix A.
The crucial aspects are: (a) the basis is stable; (b) the basis
functions have a certain order of smoothness; and (c) the basis
functions have fast decay away from an interval of length
$O(M^{-1})$ where $M$ is the number of basis functions used. Same
consistency results can be proved as long as those properties are
satisfied.

\vskip.1in\noindent{\bf Remark 2:} When $\overline{m}$, the number
of measurements is bounded, we can relax condition {\bf A2} to that
the eigenfunctions are twice continuously differentiable and with
bounded second derivative, and under this assumption we can prove a
result analogous to Theorem 1 and Corollary 1, with the
corresponding optimal rate of convergence being $n^{-2/5}$ instead
of $n^{-4/9}$.



\section{High-dimensional vector}\label{sec:matrix}

In this section, we describe a scenario where the observations are
i.i.d. Gaussian vectors, which can be approximately represented in a
known lower dimensional space (see {\bf C''}), where the effective
dimensionality of the observations grows at a rate slower than the
sample size. For convenience, we refer this setting as the
\textit{matrix case}. It can be seen that besides the proofs of the
results derived in this section sharing a lot of common features
with those in Section \ref{sec:functional}, these results also
suggest an asymptotic equivalence between the \textit{dense case}
for functional data, and the \textit{matrix case}. This means that
understanding one problem helps in understanding the other problem.
In particular, we conjecture that the results derived for the
Gaussian vectors, such as the efficient score representation
(Theorem 4), can be carried over to the functional data case with
dense measurements.

\subsection{Model}

Suppose that we have i.i.d. observations $Y_1,\cdots,Y_n$ from
$N_m(0,\overline\Sigma)$. Assume the covariance matrix
$\overline{\Sigma}$ has the following structure
$$
\overline{\Sigma} = \overline{\Sigma}_0 + \sigma^2 I_m.
$$
This may be regarded as a ``signal-plus-noise'' model, with
$\sigma^2$ representing the variance of the isotropic noise
component. We further assume that $\overline{\Sigma}_0$ has at least
$r$ positive eigenvalues, for some $r \geq 1$. The eigenvalues of
$\overline{\Sigma}_0$ are given by $s \overline{\lambda}_1 > \cdots
> s\overline{\lambda}_r > s \overline{\lambda}_{r+1} \geq \cdots
\geq 0$, where $s > 0$ is a parameter representing the ``signal
strength'' (so that $s/\sigma^2$ represents the signal-to-noise
ratio). We assume that the observations can be well represented in a
known $M$ dimensional basis $\Phi$ with $M \leq m$ (condition {\bf
C''}). Then the model space ${\cal M}_{M,r}(\Phi)$ (with $r \leq M
\leq m$) is defined as the set of all $m \times m$ matrices $\Sigma$
of the form $\Sigma = s \Phi^T B \Lambda B^T \Phi + \sigma^2 I_m$,
where $\Phi$ is an $M \times m$  matrix satisfying $\Phi \Phi^T =
I_M$, $B \in {\cal S}_{M,r}$ and $\Lambda$ is $r\times r$, diagonal
with positive diagonal elements. Note that, in order to prove
consistency of the REML estimator, we require that the intrinsic
dimension $M$ grows with $n$ sufficiently slowly. In fact, it has
been shown (e.g. in Paul (2005)) that, when $s/\sigma^2 = O(1)$, $M$
must be $o(n)$ to achieve consistency.

Throughout we assume that $\sigma^2$ and $s$ are known. Of course,
we can estimate the eigenvalues of $\overline{\Sigma}$ without any
knowledge of $s$.
The unknown parameters of the model are $B$ and $\Lambda$. The
parameter space is therefore $\Omega = {\cal S}_{M,r} \otimes
\mathbb{R}_+^r$. The estimate $(\widehat B,\widehat \Lambda)$ of
$(B,\Lambda)$ is obtained by minimizing over $\Omega$ the negative
log-likelihood (up to an additive constant and the multiplicative
factor $n$),
\begin{equation}\label{eq:restr_loglike_matrix}
L_n(B,\Lambda) =  \frac{1}{2n}\tr(\Sigma^{-1} \sum_{i=1}^n Y_i
Y_i^T) + \frac{1}{2}\log |\Sigma|.
\end{equation}
We then set the estimator of the first $r$ eigenvectors of
$\overline{\Sigma}$ as $\widehat{\Psi}=\Phi^T\widehat{B}$.

Similar to the \textit{dense case}, for asymptotic analysis, it is
more convenient to work with the transformed data $\widetilde Y_i
= \Phi Y_i$. Let $\Gamma = \Phi \Sigma \Phi^T = s B \Lambda B^T +
\sigma^2 I_M$. Then one can obtain estimates of $(B,\Lambda)$ by
minimizing over $\Omega$ the negative log-likelihood of the
transformed data:
\begin{equation}\label{eq:restr_loglike_trans}
\widetilde L_n(B,\Lambda) = \frac{1}{2n} \tr(\Gamma^{-1}
\sum_{i=1}^n \widetilde Y_i \widetilde Y_i^T) + \frac{1}{2} \log
|\Gamma|,
\end{equation}
which results in the same estimate obtained by minimizing
(\ref{eq:restr_loglike_matrix}).


\vskip.1in\noindent{\bf Remark 3:} It is known that (Muirhead,
1982), in the setting described above, the REML estimators of
$(B,\Lambda)$ coincide with the first $r$ principal components of
the sample covariance matrix of $\widetilde{Y_i} = \Phi Y_i$,
$i=1,\ldots,n$. On the other hand, based on the calculations carried
out in Appendix D, it is easy to see that the PCA estimators
$(\widehat B^{PC},\widehat \Lambda^{PC})$ satisfy the likelihood
equations $\nabla_B \widetilde L_n(\widehat B^{PC},\widehat
\Lambda^{PC}) = 0$ and $\nabla_\zeta \widetilde L_n(\widehat
B^{PC},\widehat \Lambda^{PC}) = 0$.
Thus, our approach provides an independent verification of the known
result that the PCA estimates are REML estimators under the
rank-restricted covariance model studied here.

\subsection{Consistency}\label{subsec:matrix_consistency}

We make the following assumptions about the covariance matrix.
\begin{itemize}
\item [{\bf A1'}] The eigenvalues of $\overline{\Sigma}_0$ are
given by $s\overline{\lambda}_1 \geq \cdots \geq
s\overline{\lambda}_m \geq 0$ and satisfy, for some $r \geq 1$
(fixed), (i) $\overline{c}_1 \geq \overline{\lambda}_1 > \cdots >
\overline{\lambda}_r
> \overline{\lambda}_{r+1}$ for some $\overline{c}_1 < \infty$; (ii)
$\max_{1\leq j \leq r} (\overline{\lambda}_j-
\overline{\lambda}_{j+1})^{-1} \leq \overline{c}_2 <\infty$.

\item [{\bf C''}] Assume that there exists $(B_*,\Lambda_*) \in \Omega$
(referred as ``optimal parameter'') such that, the matrix
$\Sigma_{*0} = s \Phi^T B_* \Lambda_* B_*^T \Phi$ is a close
approximation to $\overline{\Sigma}_0$ in the sense that $\beta_n :=
\parallel \overline{\Sigma}_0 - \Sigma_{*0}
\parallel_F = O(\sigma^2 \sqrt{\frac{M\log n}{n}})$.
\end{itemize}
Note that {\bf C''} implies that the observation vectors can be
closely approximated in the basis $\Phi$.

\vskip.1in\noindent {\bf Theorem 3 (matrix case):} \textit{Suppose
that {\bf A1'} and {\bf C''} hold, and $s \geq \sigma^2 > 0$.  Then
given $\eta > 0$, there exists $c_{0,\eta} > 0$ such that for
$\alpha_n = c_{0,\eta} \sigma \sqrt{\frac{M\log n}{n s}}$, with
probability at least $1-O(n^{-\eta})$, there is a minimizer
$(\widehat B,\widehat \Lambda)$ of (\ref{eq:restr_loglike_trans})
satisfying}
\begin{eqnarray*}
\parallel (I_M - B_* B_*^T) (\widehat B - B_*) \parallel_F
&\leq& \alpha_n,\\
\parallel B_*^T (\widehat B - B_*) \parallel_F &\leq&
\sqrt{\frac{s}{\sigma^2}} \alpha_n,\\
\parallel \widehat \Lambda - \Lambda_* \parallel_F &\leq&
\sqrt{\frac{s}{\sigma^2}}\alpha_n.
\end{eqnarray*}

\vskip.1in\noindent Observe that the rates obtained in Theorem 2
and Theorem 3 are identical once we replace $m$ in Theorem 2 by
$s$. Thus the number of measurements $m$ in the \textit{dense
case} is an analog of the signal strength $s$ in the
\textit{matrix case}. This important observation suggests an
\textit{asymptotic equivalence} between these two problems. This
is a result of the concentration of the matrices $\{R_i\}_{i=1}^n$
around $I_M$ for the \textit{dense case} (condition {\bf D}).
Under the matrix case, the analogs of $R_i$ exactly equal the
identity matrix. Moreover, Theorem 3 establishes the closeness of
the REML estimator to the optimal parameter, which serves as an
important step towards proving Theorem 4.

\subsection{Efficient score representation}\label{subsec:matrix_score}

When the observations are i.i.d. Gaussian vectors, we can get a more
refined result than the one stated in Theorem 3. In this section, we
show that by using the intrinsic geometry, we can get an
\textit{efficient score representation} of the REML estimator (and
hence PCA estimator). In Paul and  Johnstone (2007), a first order
approximation to the sample eigenvectors (i.e. PCA estimates) is
obtained using matrix perturbation theory (Kato, 1980).
Subsequently, it has also been shown there that the rate of
convergence of $l^2$-risk of PCA estimators is optimal. Here, we
show that the efficient score representation of the REML estimator
coincides with this first order approximation when the
signal-to-noise ratio $s/\sigma^2$ is bounded (Corollary 3). Our
approach is different from the perturbation analysis. It also
quantifies the role of intrinsic geometry of the parameter space
explicitly. Our result gives an alternative interpretation of this
approximation, and consequently, the score representation points to
an asymptotic optimality of the REML (and hence PCA) estimator.

We first introduce some notations. More details can be found in
Appendix D. Let $\zeta = \log \Lambda$ (treated interchangeably as
an $r \times 1$ vector and an $r\times r$ diagonal matrix). The
the parameter space for $(B,\zeta)$ is $\tilde{\Omega} :={\cal
S}_{M,r} \otimes \mathbb{R}^r$. Let ${\cal T}_{B} := \{U \in
\mathbb{R}^{M\times r} : B^T U = - U^T B\}$ denote the tangent
space of the \textit{Stiefel manifold} ${\cal S}_{M,r}$ at $B$.
Then the tangent space for the product manifold $\tilde{\Omega}$
at $(B,\zeta)$ is ${\cal T}_{B} \oplus \mathbb{R}^r$ (see Appendix
E for the definition of the product manifold and its tangent
space).

For notational simplicity, we use $\theta_*$ to denote
$(B_*,\zeta_*)$ and $\theta_0$ to denote $(B_0,\zeta_0)$. Define
$L(\theta_0;\theta_*) = \mathbb{E}_{\theta_*}\widetilde
L_n(\theta_0)$. Let $\nabla \widetilde L_n(\cdot)$ and $\nabla
L(\cdot; \theta_*)$ denote the \textit{intrinsic gradient} of the
functions $\widetilde L_n(\cdot)$ and $L(\cdot;\theta_*)$ with
respect to $(B,\zeta)$, respectively. Also, let $H_n(\cdot)$ and
$H(\cdot;\theta_*)$ denote the \textit{intrinsic Hessian operator}
of the functions $\widetilde L_n(\cdot)$ and $L(\cdot;\theta_*)$
with respect to $(B,\zeta)$, respectively.  Let
$H^{-1}(\cdot;\theta_*)$ denote the inverse Hessian operator of
$L(\cdot;\theta_*)$. Also we use $H_{B}(\cdot;\theta_*)$ to denote
the Hessian of $L(\cdot;\theta_*)$ w.r.t. $B$.  Notations for
Hessian w.r.t. $\zeta$ and gradients w.r.t $B$ and $\zeta$ are
defined similarly.

The following result gives the efficient score representation of the
REML estimator in the situation when $\sigma^2 = 1$ and $s=1$. The
result can be extended via rescaling to the case for arbitrary
$\sigma^2$ and $s$ with $s \geq \sigma^2 > 0$, and $s/\sigma^2$
being bounded.

\vskip.1in\noindent{\bf Theorem 4 (score representation):}
\textit{Suppose that {\bf A1'} and {\bf C''} hold with $\sigma^2 =
1$, $s=1$, and $M = o(n^{a})$ for some $a \in (0,1)$. Let $\gamma_n
= \max\{\sqrt{\frac{M \vee \log n}{n}},\beta_n\}$. Then there is a
minimizer $(\widehat B,\widehat \Lambda)$ of the negative
log-likelihood (\ref{eq:restr_loglike_trans}) such that, with
probability tending towards 1,
\begin{eqnarray}
\widehat B - B_* &=& -H_{B}^{-1}(\theta_*;\theta_*)(\nabla_B
\widetilde L_n (\theta_*)) + O(\gamma_n^2)
\label{eq:consistency_eigenvec_refined}\\
\widehat \Lambda - \Lambda_* &=& - \Lambda_*
H_{\zeta}^{-1}(\theta_*;\theta_*)(\nabla_\zeta \widetilde L_n
(\theta_*)) + O(\gamma_n^2). \label{eq:consistency_eigenval_refined}
\end{eqnarray}
In particular, from this representation, we have, with probability
tending towards 1,}
\begin{eqnarray}
\parallel \widehat B - B_* \parallel_F &=& O(\gamma_n);
\label{eq:B_hat_refined} \\
\parallel \widehat \Lambda - \Lambda_* \parallel_F &=&
O(\sqrt{\frac{\log n}{n}} + \gamma_n^2).
\label{eq:Lambda_hat_refined}
\end{eqnarray}

\vskip.1in\noindent Note that Theorem 4 gives the optimal rate of
convergence for $l^2$-risk of the estimated eigenvectors when
$\beta_n = 0$ (i.e., no model bias). This result follows from the
minimax lower bound on the risk obtained by Paul and  Johnstone
(2007). Note that, this lower bound under the current setting
follows essentially from the proof of Corollary 1. Also, when $a
\leq 1/2$, this result shows that $\widehat\Lambda$ converges at a
parametric rate. Indeed, the representation
(\ref{eq:consistency_eigenval_refined}) implies asymptotic normality
of $\widehat \Lambda$ when $a \leq 1/2$. In the derivation of
Theorem 4 we need to compute the Hessian and its inverse, which
leads to the following representation.

\vskip.1in\noindent{\bf Corollary 3:} \textit{Under the assumptions
of Theorem 4, we have the following representation:
\begin{eqnarray*}
H_{B}^{-1}(\theta_*;\theta_*)(\nabla_B \widetilde L_n(\theta_*)) &=&
[\mathbf{R}_1 \widetilde{S} B_{*1} : \cdots : \mathbf{R}_r
\widetilde{S} B_{*r}],
\end{eqnarray*}
where $B_{*j}$ is the $j$-th column of $B_*$, and
$$
\mathbf{R}_j = \sum_{1\leq i \neq j\leq r} \frac{1}{(\lambda_{*i} -
\lambda_{*j})} B_{*i} B_{*i}^T - \frac{1}{\lambda_{*j}} (I_M - B_*
B_*^T),
$$
is the resolvent operator corresponding to $\Gamma_*$ ``evaluated
at'' $(1+\lambda_{*j})$.}

\vskip.1in\noindent Combining Corollary 3 with
(\ref{eq:consistency_eigenvec_refined}), we get a first order
approximation to $\widehat B$ which coincides with the approximation
for sample eigenvectors obtained in Paul and Johnstone (2007).
However, Theorem 4 has deeper implications. Since it gives an
efficient score representation, it suggests an asymptotic optimality
of the REML estimators in the minimax sense.


\section{Proof of Theorem 1}\label{sec:proof_thm1}

Since $\sigma^2$ is fixed and assumed known, without loss of
generality, we take $\sigma^2 =1$. In this section, we give an
outline of the main ideas/steps. The details of the proofs are
given in Appendix B. The strategy of the proof is as follows. We
restrict our attention to a subset $\Theta(\alpha_n)$ of the
parameter space (referred as the \textit{restricted parameter
space}), which is the image under exponential map of the boundary
of an ellipsoid centered at 0, in the tangent space of an
``optimal parameter''. We then show that with probability tending
towards 1, for every parameter value in this restricted parameter
space, the value of the negative log-likelihood is greater than
the value of the negative log-likelihood at the optimal parameter.
Due to the Euclidean geometry of the tangent space, this implies
that with probability tending towards 1, there is a local maximum
of the log-likelihood within the image (under exponential map) of
the closed ellipsoid. The key steps of the proof are:

\begin{itemize}
\item[(i)]
Decompose the difference between the negative log-likelihood at the
optimal parameter and an arbitrary parameter in the restricted space
as a sum of three terms - a term representing the average
Kullback-Leibler divergence between the distributions, a term
representing random fluctuation in the log-likelihood, and a term
representing the \textit{model bias} (equation
(\ref{eq:loss_decomp_sparse})).

\item[(ii)]
For every fixed parameter in the restricted parameter space: (a)
provide upper and lower bounds (dependent on $\alpha_n$) for the
average Kullback-Leibler divergence; (b) provide upper bounds for
the random term and the model bias term. In both cases, the bounds
are probabilistic with exponentially small tails.

\item[(iii)]
Use a covering argument combined with a union bound to extend the
above probabilistic bounds on difference between log-likelihoods
corresponding to a single parameter in $\Theta(\alpha_n)$ to the
infimum of the difference over the entire $\Theta(\alpha_n)$.
\end{itemize}
The strategy of this proof is standard. However, in order to carry
it out we need to perform detailed computations involving the
geometry of the parameter space such as the structure of the tangent
space and the exponential map. Note that, in the current case the
geometry of the parameter space is well-understood, so that there
exist explicit form of the exponential map and a precise description
of the tangent space. This helps in obtaining the precise form of
the local Euclidean approximations around an optimal parameter in
the derivations.

\subsection{Parameter space and exponential map}

We use the following characterization of the tangent space ${\cal
T}_B$ of the Stiefel manifold ${\cal S}_{M,r}$ at a point $B$.
Any element $U \in {\cal T}_{B}$ can be expressed as $U = B A_U +
C_U$, where $A_U = - A_U^T$ and $B^T C_U = O$. We then define the
restricted parameter space centered at an optimal parameter
$(B_*,\Lambda_*)$ by
\begin{eqnarray}\label{eq:Theta_alpha_sparse}
\Theta(\alpha_n) &:=& \{(\mathbf{exp}(1,B_* A_U + C_U),\Lambda_*
\exp(D)): A_U = - A_U^T, B_*^T C_U = O, D \in \mathbb{R}^{r},
\nonumber\\
&& ~~\mbox{such that}~~~ \parallel A_U \parallel_F^2 +
\parallel C_U \parallel_F^2 +
\parallel D\parallel_F^2 = \alpha_n^2\}.
\end{eqnarray}
In the definition of $\Theta(\alpha_n)$ and henceforth, we shall
treat $\Lambda$ and $D$ interchangeably as an $r\times 1$ vector,
and an $r\times r$ diagonal matrix. The function $\mathbf{exp}(t,U)$
is the \textit{exponential map} on ${\cal S}_{M,r}$ at $B_*$,
mapping a tangent vector in ${\cal T}_{B_*}$ to a point on the
manifold.  For $U \in {\cal T}_{B_*}$ and $t \geq 0$, it is defined
as
\begin{equation*}
\mathbf{exp}(t,U) = B_* \mathbf{M}(t,U) + Q \mathbf{N}(t,U),
~~~\mbox{where}~~~
\begin{bmatrix}
\mathbf{M}(t,U)\\
\mathbf{N}(t,U)
\end{bmatrix}
= \exp \left(t\begin{bmatrix} B_*^T U & - R^T \\
R & O
\end{bmatrix} \right)
\begin{bmatrix}
I_r \\
O
\end{bmatrix},
\end{equation*}
where $\exp(\cdot)$ is the usual matrix exponential, and $QR =
(I_M - B_*B_*^T) U$ is the QR-decomposition. The properties of the
map $\mathbf{exp}(1,\cdot)$ that we shall heavily use in the
subsequent analysis (see Appendix B) are : for $U \in {\cal
T}_{B_*}$,
\begin{eqnarray}
B_*^T(\mathbf{exp}(1,U)-B_*) &=& B_*^T U + O\left(\left(\parallel
B_*^T U
\parallel_F + \parallel (I_M - B_* B_*^T)U
\parallel_F\right) \parallel U
\parallel_F\right), \label{eq:exp_U_1} \\
(I_M - B_* B_*^T)\mathbf{exp}(1,U) &=& (I_M - B_* B_*^T)U +
O\left(\parallel (I_M - B_* B_*^T)U \parallel_F \parallel U
\parallel_F\right), \label{eq:exp_U_2}
\end{eqnarray}
as $\parallel U \parallel_F \to 0$. These properties are easily
verified by using the definition of  the matrix exponential
$\exp(\cdot)$, and the Taylor series expansion.

\subsection{Loss decomposition}

We shall show that, given $\eta > 0$, for an appropriate choice of
the constant $c_{0,\eta}$ in the definition of $\alpha_n$ (in
Theorem 1), for large enough $n$, we have
\begin{equation}\label{eq:loss_bound_sparse}
\mathbb{P}\left(\inf_{(B,\Lambda) \in \Theta(\alpha_n)}
L_n(B,\Lambda)
> L_n(B_*,\Lambda_*)\right) \geq 1- O(n^{-\eta}).
\end{equation}
From this, it follows immediately that with probability tending
towards 1, there is a local minimum $(\widehat B, \widehat \Lambda)$
of $L_n(B,\Lambda)$ in the set $\overline{\Theta}(\alpha_n)$ defined
as
\begin{eqnarray*}
\overline{\Theta}(\alpha_n) &=& \{(\mathbf{exp}(1,B_* A_U +
C_U),\Lambda_* \exp(D)): A_U = - A_U^T, B_*^T C_U = O, D \in
\mathbb{R}^r
\nonumber\\
&& ~~\mbox{such that}~~~ \parallel A_U \parallel_F^2 +
\parallel C_U \parallel_F^2 +
\parallel D\parallel_F^2 \leq \alpha_n^2\},
\end{eqnarray*}
which concludes the proof of Theorem 1.

We start with the basic decomposition: {\small
\begin{eqnarray}\label{eq:loss_decomp_sparse}
&& L_n(B,\Lambda)-L_n(B_*,\Lambda_*)\nonumber\\
&=& \left[\mathbb{E}L_n(B,\Lambda) -
\mathbb{E}L_n(B_*,\Lambda_*)\right] + \left[\left(L_n(B,\Lambda) -
\mathbb{E}L_n(B,\Lambda)\right) - \left(L_n(B_*,\Lambda_*) -
\mathbb{E}L_n(B_*,\Lambda_*)\right)\right] \nonumber\\
 \nonumber\\
&=& \frac{1}{n}\sum_{i=1}^n K(\Sigma_i,\Sigma_{*i}) +
\frac{1}{2n}\sum_{i=1}^n
\tr\left((\Sigma_i^{-1}-\Sigma_{*i}^{-1})(S_i
-\overline{\Sigma}_i)\right) + \frac{1}{2n}\sum_{i=1}^n
\tr\left((\Sigma_i^{-1}-\Sigma_{*i}^{-1})(\overline{\Sigma}_i -
\Sigma_{*i})\right),
\end{eqnarray}
} where $S_i = Y_i Y_i^T$ and
\begin{eqnarray*}
K(\Sigma_i,\Sigma_{*i}) &=& \frac{1}{2}
\tr\left(\Sigma_i^{-1/2}(\Sigma_{*i} -
\Sigma_i)\Sigma_i^{-1/2}\right) - \frac{1}{2} \log|I_{m_i} +
\Sigma_i^{-1/2}(\Sigma_{*i}-\Sigma_i)\Sigma_i^{-1/2}|,
\end{eqnarray*}
is the Kullback-Leibler divergence corresponding to observation $i$.
Note that the proofs of Theorems 2 and 3 share a lot of commonality
with the \textit{sparse case} discussed here, in that these proofs
depend on the same basic decomposition of the loss function.

\subsection{Probabilistic bounds for a fixed parameter in $\Theta(\alpha_n)$}

In order to derive the results in the following three
propositions, we need to restrict our attention to an appropriate
subset of the space of the design points $\mathbf{T}$ which has
high probability. Accordingly, given $\eta
> 0$, we define such a set $A_\eta$ through
(\ref{eq:T_eta_set}) in Proposition 8 (in Appendix A). The following
proposition gives probabilistic bounds for the average
Kullback-Leibler divergence in terms of $\alpha_n$.

\vskip.1in\noindent{\bf Proposition 2:} \textit{Given $\eta > 0$,
for every $(B,\Lambda) \in \Theta(\alpha_n)$, there is a set
$A_{1,\eta}^{B,\Lambda}$ (depending on $(B,\Lambda)$), defined as
\begin{equation}\label{eq:KL_approx_2}
A_{1,\eta}^{B,\Lambda} := \left\{d_\eta^\prime \alpha_n^2 \leq
\frac{1}{n} \sum_{i=1} K(\Sigma_i,\Sigma_{*i}) \leq
d_\eta^{\prime\prime} \overline{m}^2\alpha_n^2\right\},
\end{equation}
for appropriate positive constants $d_\eta^\prime$ and
$d_\eta^{\prime\prime}$ (depending on $\overline{\lambda}_1$ and
$r$), such that for $n$ large enough, $\mathbb{P}(A_{\eta} \cap
(A_{1,\eta}^{B,\Lambda})^c) = O(n^{-(2+2\kappa)Mr -\eta})$.}

\vskip.1in\noindent Note that the bound in (\ref{eq:KL_approx_2}) is
not sharp when $\overline{m} \to \infty$, which leads the suboptimal
rates in Corollary 2. The following propositions bound the random
term and the bias term in (\ref{eq:loss_decomp_sparse}),
respectively.

\vskip.1in\noindent {\bf Proposition 3:} \textit{Given $\eta > 0$,
for each $(B,\Lambda) \in \Theta(\alpha_n)$, there is a set
$A_{2,\eta}^{B,\Lambda}$, defined as,
\begin{equation*}
A_{2,\eta}^{B,\Lambda} = \left\{\left|\frac{1}{2n}\sum_{i=1}^n
\tr((\Sigma_i^{-1}-\Sigma_{*i}^{-1})(S_i
-\overline{\Sigma}_i))\right| \leq d_{\eta} \overline{m} \alpha_n
\sqrt{\frac{M \log n}{n}} \right\},
\end{equation*}
for some $d_{\eta} > 0$, such that, $\mathbb{P}(A_{1,\eta}\cap
(A_{2,\eta}^{B,\Lambda})^c) = O(n^{-(2+2\kappa)Mr -\eta})$.}

\vskip.1in\noindent{\bf Proposition 4:} \textit{Given $\eta > 0$,
for each $(B,\Lambda) \in \Theta(\alpha_n)$, there is a set
$A_{3,\eta}^{B,\Lambda}$, defined as,
\begin{eqnarray*}
A_{3,\eta}^{B,\Lambda} = \Bigl\{ \left|\frac{1}{2n}\sum_{i=1}^n
\tr[(\Sigma_i^{-1} - \Sigma_{*i}^{-1})(\overline{\Sigma}_i -
\Sigma_{*i})]\right| &\leq& d_{\eta} \overline{m}\alpha_n
\sqrt{\frac{M \log n}{n}} \Bigr\},
\end{eqnarray*}
for some constant $d_{\eta}
> 0$, such that for large enough $n$, $\mathbb{P}(A_{\eta} \cap (A_{3,\eta}^{B,\Lambda})^c)
= O(n^{-(2+2\kappa)Mr -\eta})$.}

\vskip.1in\noindent Combining Propositions 2-4, we obtain that,
given $\eta > 0$, there is a constant $c_{0,\eta}$, such that, for
every $(B,\Lambda) \in \Theta(\alpha_n)$,
\begin{equation}\label{eq:risk_diff_prelim}
\mathbb{P}(\{L_n(B,\Lambda) - L_n(B_*,\Lambda_*) \leq \frac{1}{2}
\alpha_n^2\} \cap A_\eta) = O(n^{-(2+2\kappa)Mr -\eta}).
\end{equation}

\subsection{Covering of the space
$\Theta(\alpha_n)$}\label{subsec:covering_sparse}


To complete the proof of Theorem 1, we construct a
$\delta_n$-\textit{net} in the set $\Theta(\alpha_n)$, for some
$\delta_n
> 0$ sufficiently small. This means that, for any $(B_1,\Lambda_1)
\in \Theta(\alpha_n)$ there exists an element $(B_2,\Lambda_2)$ of
the net (with $B_k = \mathbf{exp}(1,B_* A_{U_k} + C_{U_k})$ and
$\Lambda_k = \Lambda_* \exp(D_k)$, $k=1,2$), such that we have
$\parallel B_1 - B_2 \parallel_F^2 + \parallel \Lambda_1 - \Lambda_2
\parallel_F^2 \leq \delta_n^2$.
The spaces $\{A \in \mathbb{R}^{r\times r} : A = -A^T\}$ and $\{C
\in \mathbb{R}^{M\times r} : B_*^T C = O\}$ are Euclidean
subspaces of dimension $r(r-1)/2$ and $Mr - r^2$, respectively.
Therefore, $\Theta(\alpha_n)$ is the image under
$\left(\mathbf{exp}(1,\cdot),\exp(\cdot)\right)$ of a
hyper-ellipse of dimension $p=Mr - r(r+1)/2$. Thus, using standard
construction of nets on spheres in $\mathbb{R}^p$, we can find
such a $\delta_n$-net ${\cal C}[\delta_n]$, with at most $d_1
\max\{1,(\alpha_n\delta_n^{-1})^p\}$ elements, for some $d_1 <
\infty$.

If we take $\delta_n = (\overline{m}^2 n)^{-1}$, then from
(\ref{eq:risk_diff_prelim}) using union bound it follows that, for
$n$ large enough,
\begin{equation*}
\mathbb{P}\left(\left\{\inf_{(B,\Lambda) \in {\cal
C}[\delta_n]}L_n(B,\Lambda) - L_n(B_*,\Lambda_*)
> \frac{1}{2}\alpha_n^2\right \} ~\cap~A_\eta\right) \geq 1-O(n^{-\eta}).
\end{equation*}
This result, together with the following lemma and the fact that
$\mathbb{P}(A_\eta) \geq 1 - O(n^{-\eta})$ (Proposition 8), as well
as the definition of ${\cal C}[\delta_n]$, proves
(\ref{eq:loss_bound_sparse}). The proof of Lemma 2 is given in
Appendix B.

\vskip.1in\noindent{\bf Lemma 2:} \textit{Let $(B_k,\Lambda_k)$,
$k=1,2$, be any two elements of  $\Theta(\alpha_n)$ satisfying
$\parallel B_1 - B_2\parallel_F^2 + \parallel \Lambda_1 -
\Lambda_2\parallel_F^2 \leq \delta_n^2$, with $\delta_n =
(\overline{m}^2 n)^{-1}$. Then, given $\eta > 0$, there are
constants $d_{3,\eta}, d_{4,\eta} > 0$, such that, the set
$A_{4,\eta} := \left\{\max_{1\leq i \leq n}
\parallel \overline{\Sigma}_i^{-1/2} S_i \overline{\Sigma}_i^{-1/2}
- I_{m_i} \parallel_F \leq d_{3,\eta} \overline{m} \log n
\right\}$ satisfies $\mathbb{P}(A_{4,\eta} |\mathbf{T}) \geq 1 -
O(n^{-\eta-1})$, for $\mathbf{T} \in A_{\eta}$; and on
$A_{4,\eta}$, we have $|L_n(B_1,\Lambda_1) - L_n(B_2,\Lambda_2)| =
o(\alpha_n^2)$.}

\section{Proof of Theorem 3}\label{sec:proof_thm3}

There is essentially only one step where the proof of Theorem 3
differs from that of Theorem 1. It involves providing sharper
bounds for the Kullback-Leibler divergence between an ``optimal
parameter'', and an arbitrary parameter in the restricted
parameter space $\widetilde \Theta(\alpha_n)$, an ellipsoid in the
tangent space at the ``optimal parameter'':
\begin{eqnarray*}
\widetilde \Theta(\alpha_n) &=& \{(\mathbf{exp}(1,B_* A_U +
C_U),\Lambda_* \exp(D)): A_U = - A_U^T, B_*^T C_U = O,  D \in
\mathbb{R}^r
\nonumber\\
&& ~~\mbox{such that}~~~\frac{\sigma^2}{s} \parallel
A_U\parallel_F^2 +
\parallel C_U \parallel_F^2 +
\frac{\sigma^2}{s} \parallel D\parallel_F^2 = \alpha_n^2\}.
\end{eqnarray*}
Note that now the restricted parameter space is the image (under
exponential maps) of an ellipse, whose principal axes can differ
substantially depending on the signal-to-noise ratio $s/\sigma^2$.
This is crucial for obtaining the sharper bounds for the
Kullback-Leibler divergence (see equation
(\ref{eq:norm_bound_lower_1})). As in Section \ref{sec:proof_thm1},
our strategy is to show that, given $\eta
> 0$, for an appropriate choice of $c_{0,\eta}$, for
large enough $n$, we have
\begin{equation*}
\mathbb{P}\left(\inf_{(B,\Lambda) \in \widetilde\Theta(\alpha_n)}
\widetilde L_n(B,\Lambda)
> \widetilde L_n(B_*,\Lambda_*)\right) \geq 1- O(n^{-\eta}).
\end{equation*}
From this, we conclude the proof of Theorem 3 using similar
arguments as in the proof of Theorem 1.

Define $\widetilde{S} = \frac{1}{n}\sum_{i=1}^n \widetilde Y_i
\widetilde Y_i^T$, where $\widetilde Y_i = \Phi Y_i$. Then, for an
arbitrary $(B,\Lambda) \in \widetilde \Theta(\alpha_n)$, we have
the following decomposition:
\begin{eqnarray}\label{eq:loss_decomp}
&&\widetilde L_n(B,\Lambda)-\widetilde L_n(B_*,\Lambda_*)=
\left[\left(\widetilde L_n(B,\Lambda) - \mathbb{E}\widetilde
L_n(B,\Lambda)\right) - \left(\widetilde L_n(B_*,\Lambda_*) -
\mathbb{E}\widetilde L_n(B_*,\Lambda_*)\right)\right] \nonumber\\
&=& \hskip-.1in K(\Gamma,\Gamma_{*}) + \frac{1}{2}
\tr\left((\Gamma^{-1}-\Gamma_{*}^{-1})(\widetilde{S}
-\overline{\Gamma})\right) +  \frac{1}{2}
\tr\left((\Gamma^{-1}-\Gamma_{*}^{-1})(\overline{\Gamma} -
\Gamma_{*})\right)
\end{eqnarray}
with
\begin{eqnarray*}
K(\Gamma,\Gamma_*) &=& \frac{1}{2} \tr(\Gamma^{-1}(\Gamma_* -
\Gamma)) - \frac{1}{2} \log|I_M + \Gamma^{-1}(\Gamma_*-\Gamma)| \nonumber\\
&=& \frac{1}{2} \tr(\Gamma^{-1/2}(\Gamma_* - \Gamma)\Gamma^{-1/2}) -
\frac{1}{2} \log|I_M + \Gamma^{-1/2}(\Gamma_* -
\Gamma)\Gamma^{-1/2}|,
\end{eqnarray*}
being the Kullback-Leibler divergence between the probability
distributions $N_M(0,\Gamma)$ and $N_M(0,\Gamma_*)$, where
$\Gamma^{-1/2} = (\Gamma^{1/2})^{-1}$, and $\Gamma^{1/2}$ is a
symmetric, positive definite, square root of $\Gamma$. The
following is an analogue of Proposition 2.

\vskip.1in\noindent{\bf Proposition 5:} \textit{Under the
assumptions of Theorem 3, there exist constants
$c^{\prime},c^{\prime\prime} > 0$ such that, for sufficiently large
$n$,
\begin{equation}\label{eq:Kullback_inequality_final}
c^{\prime} \alpha_n^2 \left(\frac{s}{\sigma^2}\right) \leq
K(\Gamma,\Gamma_{*}) \leq c^{\prime\prime} \alpha_n^2
\left(\frac{s}{\sigma^2}\right),
\end{equation}
for all $(B,\Lambda) \in \widetilde\Theta(\alpha_n)$, where $\Gamma
= s B \Lambda B^T + \sigma^2 I_M$ and $\Gamma_* = s B_* \Lambda_*
B_*^T + \sigma^2 I_M$.}


\vskip.1in\noindent The following are analogues of the Propositions
3 and 4, respectively.

\vskip.1in\noindent{\bf Proposition 6:} \textit{Given $\eta > 0$,
there exists a constant $c_{\eta} > 0$, such that for each
$(B,\Lambda) \in \widetilde\Theta(\alpha_n)$,
\begin{eqnarray*}
\mathbb{P}\left(\left|\tr((\Gamma^{-1}-\Gamma_{*}^{-1})(\widetilde{S}
-\overline{\Gamma}))\right| \leq c_{\eta} \sqrt{\frac{M \log
n}{n}} \sqrt{\frac{s}{\sigma^2}}\alpha_n\right) &\geq&
1-O(n^{-(2+2\kappa)Mr -\eta}).
\end{eqnarray*}
}

\vskip.1in\noindent This proposition can be easily proved using an
exponential inequality by Davidson and Szarek (2001) on the
fluctuations of the extreme eigenvalues of a Wishart matrix.

\vskip.1in\noindent{\bf Proposition 7:} \textit{There is a constant
$c > 0$ such that, uniformly over $(B,\Lambda) \in
\widetilde\Theta(\alpha_n)$,
\begin{eqnarray*}
\left|\tr((\Gamma^{-1}-\Gamma_{*}^{-1})(\overline{\Gamma} -
\Gamma_{*}))\right| &\leq& \parallel \Gamma^{-1} - \Gamma_{*}^{-1}
\parallel_F
\parallel \overline{\Gamma} -
\Gamma_{*} \parallel_F ~\leq~ c \frac{1}{\sigma^2}
\sqrt{\frac{s}{\sigma^2}} \alpha_n \beta_n.
\end{eqnarray*}
}

\vskip.1in\noindent Propositions 5-7 (together with conditions {\bf
A1'} and {\bf C''}) show that, for an appropriate choice of
$c_{0,\eta}$, $\widetilde L_n(B,\Lambda) - \widetilde
L_n(B_*,\Lambda_*) \geq c' \alpha_n^2$, for some $c' > 0$ with very
high probability, for every fixed $(B,\Lambda) \in
\widetilde\Theta(\alpha_n)$. The proof of Theorem 3 is finished by
constructing a $\delta_n$-\textit{net} similarly as in Section
\ref{subsec:covering_sparse} for the \textit{sparse case}.

\section{Proof of Theorem 4}\label{sec:proof_thm4}

The basic strategy of the proof is similar to that in classical
inference with Euclidean parameter space. The main difficulty in
present context lies in dealing with the Hessian operator of the
log-likelihood (intrinsic Fisher information operator) and its
inverse. Details of these calculations are given in Appendix D.

Rewrite the negative log-likelihood (\ref{eq:restr_loglike_trans})
(up to a multiplicative constant) as
\begin{equation}\label{eq:loss_vector}
\widetilde L_n(B,\Lambda) = \tr(\Gamma^{-1}\widetilde{S}) + \log
|\Gamma|, ~~~\mbox{where}~~\widetilde{S} = \frac{1}{n}\sum_{i=1}^n
\widetilde Y_i \widetilde Y_i^T.
\end{equation}
By Theorem 3, given $\eta > 0$, there is a constant $c_{3,\eta}
>0$ such that the set
$$
\widetilde A_{3,\eta} := \left\{ \parallel \widehat U \parallel_F^2
+ \parallel \widehat D \parallel_F^2 \leq
c_{3,\eta}\alpha_n^2\right\}
$$
has probability at least $1- O(n^{-\eta})$,  where $\alpha_n =
c_{0,\eta} \sqrt{\frac{M\log n}{n}}$, and $(\widehat U, \widehat
D) \in {\cal T}_{B_*} \oplus \mathbb{R}^r$ is such that $(\widehat
B, \widehat \Lambda) := (\mathbf{exp}(1,\widehat
U),\Lambda_*\exp(\widehat D))$ is a minimizer of
(\ref{eq:loss_vector}).

First, by the same concentration bound for singular values of random
matrices with i.i.d. Gaussian entries (Davidson and Szarek, 2001)
used in the proof of Proposition 6, there exists $c_{4,\eta} > 0$,
such that the set
$$\widetilde
A_{4,\eta} := \left\{\parallel \widetilde{S} -
\overline{\Gamma}\parallel \leq c_{4,\eta}\sqrt{\frac{M \vee \log
n}{n}}\right\}$$ has probability at least $1 -O(n^{-\eta})$. It
then follows that, we can choose an appropriate constant
$c_{5,\eta}
>0$ such that, on $\widetilde A_{3,\eta} \cap \widetilde
A_{4,\eta}$, $\parallel \nabla \widetilde L_n(\theta_*)\parallel
\leq c_{5,\eta} \gamma_n$, where $\gamma_n = \max\{\sqrt{\frac{M
\vee \log n}{n}},\beta_n\}$ and $\theta_*=(B_*,\Lambda_*)$. Next,
for any $X = (X_B,X_{\zeta}) \in {\cal T}_{B_*} \oplus
\mathbb{R}^r$, define
 $$\parallel X \parallel := \Bigl[\parallel
X_B\parallel_F^2 +
\parallel X_{\zeta} \parallel_F^2 \Bigr]^{1/2}.$$
Also, let $\langle \cdot,\cdot \rangle_g$ denote the canonical
metric on ${\cal T}_{B_*} \oplus \mathbb{R}^r$ (see Appendix D).
Using the fact that $\nabla \widetilde L_n(\widehat \theta) = 0$,
where $\widehat \theta = (\widehat B,\widehat \Lambda)$, and
defining $\widehat \Delta := (\widehat U, \widehat D)$, then on
$\widetilde A_{3,\eta} \cap \widetilde A_{4,\eta}$, for any $X \in
{\cal T}_{B_*} \oplus \mathbb{R}^r$ with $\parallel X
\parallel \leq 1$,
\begin{eqnarray}\label{eq:gradient_expand}
-\langle \nabla \widetilde L_n(\theta_*), X \rangle_g &=& \langle
\nabla \widetilde L_n(\widehat \theta) - \nabla \widetilde
L_n(\theta_*),
X \rangle_g \nonumber\\
&=& \langle H_n(\theta_*)(\widehat \Delta), X\rangle_g +
O(\parallel \widehat \Delta\parallel^2) + O(\gamma_n
\parallel \widehat \Delta\parallel)
\nonumber\\
&=& \langle H(\theta_*;\theta_*)(\widehat \Delta), X \rangle_g +
\langle [H_n(\theta_*) - H(\theta_*;\theta_*)](\widehat \Delta), X
\rangle_g + O(\alpha_n^2 + \alpha_n\gamma_n),
\end{eqnarray}
where $H_n(\cdot)(\widehat \Delta)$ and $H(\cdot;\theta_*)(\widehat
\Delta)$ are the corresponding \textit{covariant derivatives} of
$\widetilde L_n(\cdot)$ and $L(\cdot;\theta_*)$ in the direction of
$\widehat \Delta$.  By simple calculations based on the expressions
in Appendix D, there exists a constant $c_{6,\eta} > 0$, such that
on $\widetilde A_{3,\eta}\cap \widetilde A_{4,\eta}$, $\parallel
H_n(\theta_*)(\widehat \Delta)-H(\theta_*;\theta_*)(\widehat
\Delta)\parallel \leq c_{6,\eta} \alpha_n \gamma_n$. It can be
checked using assumptions {\bf A1'} and {\bf C''}  that the linear
operator $H^{-1}(\theta_*;\theta_*) : {\cal T}_{B_*} \oplus
\mathbb{R}^r \to {\cal T}_{B_*}\oplus \mathbb{R}^r$, is bounded in
operator norm (Appendix D). Therefore, using the definition of
covariant derivative and inverse of Hessian, from
(\ref{eq:gradient_expand}) we have, on $\widetilde A_{3,\eta}\cap
\widetilde A_{4,\eta}$,
\begin{equation}\label{eq:delta_hat_expand}
\widehat \Delta = - H^{-1}(\theta_*;\theta_*) (\nabla \widetilde
L_n(\theta_*)) + O(\alpha_n\gamma_n) + O(\alpha_n^2).
\end{equation}
Hence,  on $\widetilde A_{3,\eta}\cap \widetilde A_{4,\eta}$, the
bound on $\parallel \widehat \Delta\parallel$ can be improved from
$O(\alpha_n)$ to
\begin{equation}\label{eq:delta_hat_rate}
\parallel \widehat\Delta \parallel = O(\gamma_n\alpha_n + \alpha_n^2).
\end{equation}
We can then repeat exactly the same argument, by using
(\ref{eq:gradient_expand}) to derive (\ref{eq:delta_hat_expand}),
but now with the bound on $\parallel \widehat\Delta
\parallel$ given by (\ref{eq:delta_hat_rate}).
Since $M = O(n^{a})$ for some $a < 1$, so that $\alpha_n^2 =
o(\gamma_n)$, this way we get the more precise expression,
\begin{equation}\label{eq:delta_hat_expand_refined}
\widehat \Delta = - H^{-1}(\theta_*;\theta_*) (\nabla \widetilde
L_n(\theta_*)) + O(\gamma_n^2).
\end{equation}
Moreover, it can be easily verified that,
\begin{equation*}
\frac{\partial}{\partial \zeta} \nabla_B L(\theta_*;\theta_*) :=
\mathbb{E}_{\theta_*} \left[\frac{\partial}{\partial \zeta} \nabla_B
\widetilde L_n(\theta_*)\right]= 0.
\end{equation*}
Hence, by (\ref{eq:prod_Hessian_operator}) in Appendix E, the
Hessian operator, and its inverse, are ``block diagonal'', on the
parameter space (viewed as a product manifold), with diagonal
blocks corresponding to Hessians (inverse Hessians) w.r.t. $B$ and
$\zeta$, respectively. This yields
(\ref{eq:consistency_eigenvec_refined}) and
(\ref{eq:consistency_eigenval_refined}) in Theorem 4. Also,
(\ref{eq:B_hat_refined}) and (\ref{eq:Lambda_hat_refined}) follow
immediately from (\ref{eq:delta_hat_expand_refined}).


\section{Discussion}\label{sec:discussion}

In this paper, we have demonstrated the effectiveness of utilizing
the geometry of the non-Euclidean parameter space in determining
consistency and rates of convergence of the REML estimators of
principal components. We first study the REML estimators of
eigenvalues and eigenfunctions of the covariance kernel for
functional data, estimated from sparse, irregular measurements. The
convergence rate of the estimated eigenfunctions is shown to be
near-optimal when the number of measurements per curve is bounded
and when $M$, the number of basis functions, varies with $n$ at an
appropriate rate (Theorem 1 and Corollary 1). The technique used in
proving Theorem 1 is most suitable for dealing with the very sparse
case (i.e., number of measurements per curve is bounded). We have
also used it to prove consistency for the case where the number of
measurements increases slowly with sample size (Corollary 2).
However, this does not result in the optimal convergence rate. The
latter case is more difficult because of the complications of
dealing with inverses of random matrices ($\Sigma_i$) of growing
dimensions. A more delicate analysis, that can handle this issue
more efficiently, is likely to give tighter bounds for the average
Kullback-Leibler divergence than that obtained in Proposition 2.
Then it may be possible to extend the current technique to prove
optimality of the REML estimators in a broader regime. A variant of
the technique used for proving Theorem 1 also gives consistency of
the REML estimator for functional data in a regime of dense
measurements, as well as for a class of high-dimensional Gaussian
vectors (Theorems 2 and 3). In the latter case, we also derive an
efficient score representation (Theorem 4), which involves
determining the intrinsic Fisher information operator and its
inverse.



Now we present some conjectures we aim to pursue. First, as
discussed earlier, based on the score representation, we conjecture
the asymptotic optimality of the REML estimator for the
\textit{matrix case}. Secondly, we conjecture that there exists an
efficient score representation of the REML estimator in the
functional data problem as well. If so, then this estimator is
likely to achieve the optimal nonparametric rate (for a broader
regime), and may even be asymptotically optimal. This may explain
the superior numerical performance of the REML estimator observed by
Peng and Paul (2007). Thirdly, our results (Theorems 2 and 3) give a
strong indication of an asymptotic equivalence between two classes
of problems : statistical inference for functional data with dense
measurements; and inference for high-dimensional i.i.d. Gaussian
vectors. Finally, in this paper we have not addressed the issue of
model selection. A procedure for selection of $M$ and $r$, based on
an approximate leave-one-curve-out cross-validation score, has been
proposed and implemented in Peng and Paul (2007). This approximation
is based on a second order Taylor expansion of the negative
log-likelihood at the estimator and it involves the intrinsic Fisher
information operator and its inverse. Therefore, based on the
analysis presented here, it is conjectured that the approximate CV
score thus defined is asymptotically consistent for the class of
models considered in this paper.

\section*{Appendix A : Properties of cubic B-spline basis}

In many proofs of this paper, we need to use some properties of the
cubic $B$-spline basis. We state some of them. More details can be
found in de Boor (1978) and deVore and Lorentz (1993). Let
$\tilde{\boldsymbol{\phi}} = (\tilde \phi_1,\ldots,\tilde \phi_M)^T$
be the (standard) cubic B-spline basis functions on $[0,1]$ with
equally spaced knots. Then, the orthonormalized spline functions
$\phi_1,\ldots,\phi_M$ are defined through $\boldsymbol{\phi}(t) =
G_{\phi,M}^{-1/2} \tilde{\boldsymbol{\phi}}(t)$, where $G_{\phi,M}
:= ((\int \tilde\phi_k(t)\tilde \phi_l(t) dt))_{k,l=1}^M$, is the
\textit{Gram matrix} of $\tilde{\boldsymbol{\phi}}$. It is known
(cf. de Boor (1974), Burman (1985)) that $G_{\phi,M}$  is an
$M\times M$ banded matrix, and satisfies,
\begin{equation}\label{eq:Gram_spline}
\frac{c_{\phi,0}}{M} I_M \leq  G_{\phi,M} \leq \frac{c_{\phi,1}}{M}
I_M~~~\mbox{for some constants}~~0 < c_{\phi,0} < c_{\phi,1} <
\infty.
\end{equation}
From this, and other properties of cubic B-splines (deVore and
Lorentz, 1993, Chapter 13), we also have the following:
\begin{itemize}
\item[{\bf S1}] $\sup_{t \in [0,1]} \sum_{k=1}^M \phi_k^2(t) \leq c_{\phi,2}
M$ for some constant $c_{\phi,2} > 0$.

\item[{\bf S2}] For any function $f \in C^{(4)}([0,1])$, we
have $\parallel f - P_{\phi,M}(f)\parallel_\infty = \parallel
f^{(4)}\parallel_\infty O(M^{-4})$, where $P_{\phi,M}(f) =
\sum_{k=1}^M \langle f,\phi_k\rangle \phi_k$ denotes the projection
of $f$ onto span$\{\phi_1,\ldots,\phi_M\}=$
span$\{\tilde\phi_1,\ldots,\tilde\phi_M\}$.
\end{itemize}
Note that, property {\bf S2} and assumption {\bf A2} imply the
existence of orthonormal functions $\{\psi_{*k}\}_{k=1}^r$ of the
from
$$
(\psi_{*1}(t),\ldots,\psi_{*r}(t)) = (\boldsymbol{\psi}_*(t))^T =
B_*^T \boldsymbol{\phi}(t),~~B_*^T B_* = I_r,
$$
which  satisfy
\begin{equation}\label{eq:psi_spline_approx}
\max_{1\leq k \leq r} \parallel \overline{\psi}_k -
\psi_{*k}\parallel_\infty \leq c_{\phi,3} M^{-4} \max_{1\leq k \leq
r}\parallel \overline{\psi}_k^{(4)} \parallel_\infty.
\end{equation}
Using these properties we obtain the following approximation to the
important quantity $\parallel \Phi_i\parallel$, where $\Phi_i =
[\boldsymbol{\phi}(T_{i1}):\ldots:\boldsymbol{\phi}(T_{im_i})]$ and
$\parallel \cdot \parallel$ denotes the operator norm. This result
will be extremely useful in the subsequent analysis.

\vskip.1in\noindent{\bf Proposition 8:} \textit{Given $\eta > 0$,
there is an event $A_{\eta}$ defined in terms of the design points
$\mathbf{T}$, with probability at least $1- O(n^{-\eta})$, such
that on the set $A_{\eta}$,
\begin{eqnarray}\label{eq:Phi_random}
\parallel \Phi_i \parallel^2 &\leq&  \overline{m}c_{g,1} + \sqrt{5}
c_{\phi,0}^{-1} d_{\eta} [(M^{3/2}\log n) \vee (M
\sqrt{\overline{m}\log n})],
\end{eqnarray}
for some constant $d_\eta > 0$, and $c_{g,1}$ is the constant in
condition {\bf B2}. Furthermore, for all $\mathbf{T}$, we have the
non-random bound}
\begin{eqnarray}\label{eq:Phi_nonrandom}
\parallel \Phi_i \parallel^2 &\leq& c_{\phi,2} \overline{m} M,
~~~\mbox{for all}~~~i=1,\ldots,n.
\end{eqnarray}

\vskip.1in\noindent{\bf Proof :} First, (\ref{eq:Phi_nonrandom})
follows from the bound {\bf S1}, since
\begin{eqnarray*}
\parallel \Phi_i \parallel^2 &=& \parallel \Phi_i^T \Phi_i \parallel
= \parallel \Phi_i \Phi_i^T \parallel \leq \parallel
\Phi_i\parallel_F^2 = \tr (\Phi_i^T \Phi_i) \nonumber\\
&=& \sum_{j=1}^{m_i} \sum_{k=1}^M (\phi_k(T_{ij}))^2 \leq c_{\phi,2}
m_i M \leq c_{\phi,2} \overline{m} M,~~\mbox{for all}~i=1,\ldots,n.
\end{eqnarray*}
In order to prove (\ref{eq:Phi_random}), first write
\begin{equation}\label{eq:Phi_diff_repr}
\frac{1}{m_i} \Phi_i \Phi_i^T  - \int
\boldsymbol{\phi}(t)(\boldsymbol{\phi}(t))^T g(t) dt =
G_{\phi,M}^{-1/2} \left[((\frac{1}{m_i}\sum_{j=1}^{m_i}
[\tilde\phi_k(T_{ij}) \tilde\phi_l(T_{ij}) -
\mathbb{E}(\tilde\phi_k(T_{ij}) \tilde\phi_l(T_{ij}))]))_{k,l=1}^M
\right] G_{\phi,M}^{-1/2}.
\end{equation}
Next, observe that, $\mathbb{E}[\phi_k(T_{i1})\phi_l(T_{i1})]^2 =
\int (\tilde\phi_k(t))^2 (\tilde\phi_l(t))^2 g(t) dt=0$  for
$|k-l|
> 3$; and is within $[c_{\phi,4} c_{g,0} M^{-1}, c_{\phi,5} c_{g,1}
M^{-1}]$, for constants $0 < c_{\phi,4} < c_{\phi,5} < \infty$, if
$|k-l| \leq 3$. Then, using the fact that $\max_{1\leq k \leq M}
\parallel \tilde \phi_k
\parallel_\infty$ is bounded,
it follows from Bernstein's inequality that the set $A_{\eta}$
defined by
\begin{equation}\label{eq:T_eta_set}
A_{\eta} = \left\{\mathbf{T} : \max_{1\leq i \leq n} \max_{1\leq
k,l \leq M} \left|\frac{1}{m_i}\sum_{j=1}^{m_i}
[\tilde\phi_k(T_{ij}) \tilde\phi_l(T_{ij}) -
\mathbb{E}(\tilde\phi_k(T_{ij}) \tilde\phi_l(T_{ij}))]\right| \leq
d_{1,\eta} \left(\frac{\log n}{\underline{m}}\right) \vee
\sqrt{\frac{\log n}{\underline{m} M}}\right\}
\end{equation}
has probability at least $1-O(n^{-\eta})$, for some constant
$d_{1,\eta} > 0$. Now, we can bound the Frobenius norm of the
matrix in  (\ref{eq:Phi_diff_repr}) by using
(\ref{eq:Gram_spline}) and (\ref{eq:T_eta_set}), and the fact that
the matrix has $O(M)$ nonzero elements. Then using
(\ref{eq:g_cond}) we derive (\ref{eq:Phi_random}).

\section*{Appendix B : Proofs for the sparse case}

{\bf Proof of Proposition 2 :}  The main challenge in the proof of
Proposition 2 is to efficiently approximate the average
Kullback-Leibler divergence. We can express
$K(\Sigma_i,\Sigma_{*i})$ as
\begin{equation}\label{eq:Kullback_expr}
K(\Sigma_i,\Sigma_{*i}) = \frac{1}{2} \sum_{j=1}^m
[\lambda_j(R_{*i}) - \log(1+\lambda_j(R_{*i}))],
\end{equation}
where $\lambda_j(R_{*i})$ is the $j$-th largest eigenvalue of
$R_{*i} = \Sigma_i^{-1/2}(\Sigma_{*i} - \Sigma_i)\Sigma_{i}^{-1/2}$.
Using the inequality $e^x \geq 1+x$ for $x \in \mathbb{R}$ (so that
each term in the summation in (\ref{eq:Kullback_expr}) is
nonnegative), and the Taylor series expansion for $\log(1+x)$ for
$|x| < 1$, it can be shown that, given $\epsilon > 0$ sufficiently
small (but fixed), there exist constants $0 < c_{1,\epsilon} <
c_{2,\epsilon} < \infty$ such that for $\parallel R_{*i}
\parallel_F \leq \epsilon$,
\begin{equation}\label{eq:Kullback_inequality_sparse}
c_{1,\epsilon} \parallel R_{*i}
\parallel_F^2 \leq K(\Sigma_i,\Sigma_{*i}) \leq c_{2,\epsilon}\parallel
R_{*i}
\parallel_F^2.
\end{equation}
Next, observe that
\begin{equation*}
\frac{\parallel \Sigma_{*i}^{-1/2}(\Sigma_i
-\Sigma_{*i})\Sigma_{*i}^{-1/2}
\parallel_F}{1+\parallel \Sigma_{*i}^{-1/2}(\Sigma_i
-\Sigma_{*i})\Sigma_{*i}^{-1/2} \parallel_F} \leq ~
\parallel R_{*i} \parallel_F ~ \leq   \frac{\parallel \Sigma_{*i}^{-1/2}(\Sigma_i
-\Sigma_{*i})\Sigma_{*i}^{-1/2}\parallel_F}{1-\parallel
\Sigma_{*i}^{-1/2}(\Sigma_i -\Sigma_{*i})\Sigma_{*i}^{-1/2}
\parallel_F},
\end{equation*}
whenever $\parallel \Sigma_{*i}^{-1/2}(\Sigma_i
-\Sigma_{*i})\Sigma_{*i}^{-1/2} \parallel_F < 1$.  the proof of
Proposition 1 can thus be reduced to finding probabilistic bounds
for $\frac{1}{n}\sum_{i=1}^n \parallel \Sigma_{*i}^{-1/2}(\Sigma_i -
\Sigma_{*i}) \Sigma_{*i}^{-1/2}\parallel_F^2$.

One difficulty in obtaining those bounds is in handling the inverse
of the matrices $\Sigma_{*i}$. In order to address that, and some
related issues, we use the properties of the cubic spline basis
derived in Appendix A. In the following lemmas we confine ourselves
to the restricted parameter space $\Theta(\alpha_n)$, i.e.,
$(B,\Lambda) \in \Theta(\alpha_n)$.

\vskip.1in\noindent{\bf Lemma 3:} \textit{Under the assumptions of
Theorem 1 (for $m_i$'s bounded), or Corollary 2 (for $m_i$'s
increasing slowly with $n$),
\begin{equation}\label{eq:Sigma_quad_bound}
(1 + d_1 \overline{\lambda}_1 r \overline{m})^{-1}
\parallel \Sigma_i -
\Sigma_{*i} \parallel_F \leq \parallel \Sigma_{*i}^{-1/2}(\Sigma_i -
\Sigma_{*i})\Sigma_{*i}^{-1/2}\parallel_F \leq \parallel \Sigma_i -
\Sigma_{*i} \parallel_F,~~\mbox{for all}~i=1,\ldots,n,
\end{equation}
for some constant $d_1 > 0$}

\vskip.1in\noindent{\bf Proof :} From condition {\bf A2} and
(\ref{eq:psi_spline_approx}), it follows that, $\exists~D_1 > 0$
such that, for all $M$,
\begin{equation}\label{eq:psi_star_bound}
\max_{1\leq k \leq r} \parallel \psi_{*k}\parallel_\infty \leq D_1 <
\infty.
\end{equation}
This, together with the definition of $(B_*,\Lambda_*)$, leads to
the following bound on the eigenvalues of the matrices
$\Sigma_{*i}$:
\begin{equation}\label{eq:Sigma_eigen_bound}
1 \leq \lambda_{\min}(\Sigma_{*i}) \leq \lambda_{\max}(\Sigma_{*i})
\leq 1 + D_1 r m_i  \lambda_{*1}  \leq 1 + d_1 \overline{\lambda}_1
r \overline{m}, ~~~\mbox{for all}~i=1,\ldots,n,
\end{equation}
for some $d_1 > 0$, from which (\ref{eq:Sigma_quad_bound}) follows.

\vskip.1in\noindent{\bf Lemma 4:} \textit{Under the assumptions of
Theorem 1 (for $m_i$'s bounded), or Corollary 2 (for $m_i$'s
increasing slowly with $n$), given any $\eta > 0$, on the event
$A_{\eta}$  defined through (\ref{eq:T_eta_set}) in Proposition 8,
which has probability at least $1-O(n^{-\eta})$, for sufficiently
large $n$,
\begin{eqnarray}\label{eq:Sigma_diff_upper2}
\max_{1\leq i \leq n} \parallel \Sigma_i -
\Sigma_{*i}\parallel_F^2 &\hskip-.1in \leq & \hskip-.1in
\left[d_{3,\eta}  \left(1+ d_{1}\left[(\frac{M^{3/2}\log
n}{\overline{m}}) \vee \sqrt{\frac{M^2\log
n}{\overline{m}}}\right]\right) \overline{m}^2 \alpha_n^2\right]
\wedge \left[ d_{2} M \overline{m}^2 \alpha_n^2\right],
\end{eqnarray}
where the second bound holds for all $\mathbf{T}$. Here
$d_1,d_2,d_{3,\eta} >0$ are appropriate constants depending on $r$
and $\overline{\lambda}_1$.}

\vskip.1in\noindent{\bf Proof :} An upper bound for $\parallel
\Sigma_i - \Sigma_{*i}\parallel_F$ is obtained by expressing $I_M =
B_* B_*^T + (I_M - B_* B_*^T)$, and then applying the triangle
inequality,
\begin{eqnarray}\label{eq:Sigma_diff_upper1}
&& \parallel \Sigma_i - \Sigma_{*i}\parallel_F  = \parallel \Phi_i^T
(B\Lambda B^T - B_*
\Lambda_* B_*^T) \Phi_i\parallel_F \nonumber\\
&\leq&
\parallel \Phi_i^T B_* (B_*^T B \Lambda B^T B_*  - \Lambda_* )B_*^T \Phi_i
\parallel_F + 2 \parallel \Phi_i^T B_* B_*^T B \Lambda B^T (I_M - B_* B_*^T)
\Phi_i\parallel_F \nonumber\\
&&  + \parallel \Phi_i^T (I_M - B_* B_*^T) B \Lambda B^T (I_M - B_*
B_*^T) \Phi_i \parallel_F \nonumber\\
&\leq& \parallel \Phi_i^T B_* \parallel^2 \parallel  B_*^T B \Lambda
B^T B_* - \Lambda_* \parallel_F + 2 \parallel \Phi_i^T B_* \parallel
\parallel \Lambda \parallel \parallel B^T (I_M - B_* B_*^T) \Phi
\parallel_F \nonumber\\
&& + \parallel \Lambda \parallel
\parallel \Phi_i^T (I_M - B_* B_*^T) B \parallel_F^2 \nonumber\\
&\leq& D_1 r \overline{m} \parallel  B_*^T B \Lambda B^T B_* -
\Lambda_* \parallel_F \nonumber\\
&& + \sqrt{d_4 r \overline{m}} \overline{\lambda}_1
\parallel \Phi_i^T (I_M - B_* B_*^T) B \parallel_F \left[1+ (D_1 r
\overline{m})^{-1/2} \parallel \Phi_i^T (I_M - B_* B_*^T) B
\parallel_F\right],
\end{eqnarray}
for some $d_4 > 1$. For the second inequality we use $\parallel
B_*^T B \parallel \leq 1$, and for the last inequality we use
(\ref{eq:psi_star_bound}) and (\ref{eq:Sigma_eigen_bound}). Next, by
using (\ref{eq:Phi_nonrandom}), (\ref{eq:exp_U_2}) and
(\ref{eq:M_condition_sparse}), we obtain the (nonrandom) bound
\begin{equation}\label{eq:ortho_proj_B_bound}
\overline{m}^{-1} \max_{1\leq i \leq n} \parallel \Phi_i^T (I_M -
B_* B_*^T) B
\parallel_F^2 \leq c_{\phi,2} M  \parallel (I_M - B_*
B_*^T) B
\parallel_F^2 \leq  c_{\phi,2} M \alpha_n^2
(1+o(1)) = o(1).
\end{equation}
Then the bound in (\ref{eq:Sigma_diff_upper1}) can be majorized by,
$$
D_1 r \overline{m} \parallel  B_*^T B \Lambda B^T B_* - \Lambda_*
\parallel_F  + \sqrt{d_4 r
\overline{m}} \overline{\lambda}_1
\parallel \Phi_i \parallel \parallel (I_M - B_* B_*^T) B \parallel_F(1+o(1)).
$$
Using (\ref{eq:Phi_random}) to bound $\parallel \Phi_i \parallel$,
from (\ref{eq:Sigma_diff_upper1}), and the definition of
$\Theta(\alpha_n)$ together with (\ref{eq:exp_U_1}) and
(\ref{eq:exp_U_2}), we obtain (\ref{eq:Sigma_diff_upper2}).

\vskip.1in\noindent{\bf Lemma 5:} \textit{Under the assumptions of
Theorem 1 (for $m_i$'s bounded), or Corollary 2 (for $m_i$'s
increasing slowly with $n$), for any given $\eta >0$, there is a
positive sequence $\varepsilon_{1,n} = o(1)$ (depending on $\eta$),
and a constant $d_{1,\eta}
> 0$, such that
\begin{equation*}
\mathbb{P}\left[\frac{1}{n}\sum_{i=1}^n \parallel \Sigma_i -
\Sigma_{*i}\parallel_F^2
> d_{1,\eta} \underline{m}^2 \alpha_n^2
(1-\varepsilon_{1,n})\right] \geq 1 - O(n^{-(2+2\kappa) Mr -
\eta}),
\end{equation*}
where $\kappa$ is as in {\bf B1}.}

\vskip.1in\noindent{\bf Proof :} Let $\Delta = B \Lambda B^T - B_*
\Lambda_* B_*^T$. Observe that, by definition of $\Theta(\alpha_n)$
(equation (\ref{eq:Theta_alpha_sparse})), and equations
(\ref{eq:exp_U_1}) and (\ref{eq:exp_U_2}), for large $n$,
\begin{equation}\label{eq:Delta_alpha_rate}
\parallel \Delta \parallel_F^2 \leq c_* \alpha_n^2, ~~~~~(\mbox{for some
constant} ~~c_* > 0).
\end{equation}
First, consider the lower bound
\begin{eqnarray}\label{eq:Sigma_diff_lower1}
\parallel \Sigma_i -
\Sigma_{*i}\parallel_F^2 &=& \tr[\Phi_i^T \Delta  \Phi_i \Phi_i^T
\Delta \Phi_i] ~\geq~ \sum_{j_1\neq j_2}^{m_i}
[(\boldsymbol{\phi}(T_{ij_1}))^T \Delta
\boldsymbol{\phi}(T_{ij_2})]^2.
\end{eqnarray}
We are going to derive an exponential tail bound for $\frac{1}{n}
\sum_{i=1}^n \sum_{j_1\neq j_2}^{m_i}
[(\boldsymbol{\phi}(T_{ij_1}))^T \Delta
\boldsymbol{\phi}(T_{ij_2})]^2$. Rewriting the term on the extreme
right, and using the fact that $\{T_{ij}\}_{j=1}^{m_i}$ are i.i.d.
with density $g : c_{g,0} \leq g \leq c_{g,1}$ ({\bf B2}), we have,
for all $i$,
\begin{eqnarray}\label{eq:Sigma_diff_lower2}
\mathbb{E} \sum_{j_1\neq j_2}^{m_i} [(\boldsymbol{\phi}(T_{ij_1}))^T
\Delta \boldsymbol{\phi}(T_{ij_2})]^2 &=& m_i(m_i-1)
\tr\left(\mathbb{E}[\boldsymbol{\phi}(T_{i1})(\boldsymbol{\phi}(T_{i1}))^T]
\Delta
\mathbb{E}[\boldsymbol{\phi}(T_{i2})(\boldsymbol{\phi}(T_{i2}))^T]
\Delta\right)\nonumber\\
&\in& \left(c_{g,0}^2 \underline{m}(\underline{m}-1)
\parallel \Delta
\parallel_F^2 ,~~ c_{g,1}^2 \overline{m}(\overline{m}-1)
\parallel \Delta
\parallel_F^2\right)  \nonumber\\
&\in& \left(d_1^\prime \underline{m}^2 \alpha_n^2 (1+o(1)),~~~
d_1^{\prime\prime} \overline{m}^2 \alpha_n^2 (1+o(1))\right),
\end{eqnarray}
for some $d_1^{\prime\prime} \geq d_1^\prime > 0$ (whose values
depend on $c_{g,0}$, $c_{g,1}$ and the constants appearing in {\bf
A1}), where in the last step we use (\ref{eq:Delta_alpha_rate}). The
last inequality uses (\ref{eq:weilandt}),
(\ref{eq:eigenvec_perturb}), the definition of $\Theta(\alpha_n)$,
and the properties (\ref{eq:exp_U_1}) and (\ref{eq:exp_U_2}). Notice
that, the variance of $\sum_{j_1\neq j_2}^{m_i}
[(\boldsymbol{\phi}(T_{ij_1}))^T \Delta
\boldsymbol{\phi}(T_{ij_2})]^2$ can be bounded, for sufficiently
large $n$, as
\begin{eqnarray*}
&& \max_{1\leq i \leq n} \mbox{Var}\left(\sum_{j_1\neq j_2}^{m_i}
[(\boldsymbol{\phi}(T_{ij_1}))^T \Delta
\boldsymbol{\phi}(T_{ij_2})]^2\right) \nonumber\\
&\leq& \max_{1\leq i \leq n} \mathbb{E}\left(\parallel \Sigma_i -
\Sigma_{*i}\parallel_F^2 \sum_{j_1\neq j_2}^{m_i}
[(\boldsymbol{\phi}(T_{ij_1}))^T \Delta
\boldsymbol{\phi}(T_{ij_2})]^2\right)
\nonumber\\
&\leq& d_2' M (\overline{m}^2 \alpha_n^2)^2 ~=:~ V_{1,n},
\end{eqnarray*}
where $d_2' > 0$ is some constant. In the above, we obtain the
first inequality by (\ref{eq:Sigma_diff_lower1}), and the second
inequality by using (\ref{eq:Sigma_diff_upper2}) and
(\ref{eq:Sigma_diff_lower2}). Next, $\sum_{j_1\neq j_2}^{m_i}
[(\boldsymbol{\phi}(T_{ij_1}))^T \Delta
\boldsymbol{\phi}(T_{ij_2})]^2$, for $i=1,\ldots,n$, are
independent, and bounded  by $K_{1,n} := d_4 M \overline{m}^2
\alpha_n^2$, for a constant $d_4 > 0$ (using
(\ref{eq:Sigma_diff_upper2})). Hence, by applying Bernstein's
inequality, and noticing that $K_{1,n} \sqrt{\frac{M \log n}{n}} =
o(\sqrt{V_{1,n}})$, and $\sqrt{V_{1,n}}\sqrt{\frac{M \log n}{n}} =
o(\overline{m}^2 \alpha_n^2)$  (by (\ref{eq:M_condition_sparse}),
or (\ref{eq:m_M_condition})), the result follows.

\vskip.1in\noindent{\bf Lemma 6:} \textit{Under the assumptions of
Theorem 1 (for $m_i$'s bounded), or Corollary 2 (for $m_i$'s
increasing slowly with $n$), for any given $\eta >0$, there is a
positive sequence $\varepsilon_{2,n} = o(1)$ and a constant
$d_{2,\eta} > 0$, such that
\begin{equation}\label{eq:Sigma_diff_upper_concen}
\mathbb{P}\left[\frac{1}{n}\sum_{i=1}^n \parallel \Sigma_i -
\Sigma_{*i}\parallel_F^2 < d_{2,\eta} \overline{m}^2 \alpha_n^2
(1+\varepsilon_{2,n})\right] \geq 1 - O(n^{-(2+2\kappa)Mr -
\eta}),
\end{equation}
where $\kappa$ is as in {\bf B1}.}

\vskip.1in\noindent{\bf Proof :} From the proof of Lemma 4,
especially the inequalities (\ref{eq:Sigma_diff_upper1}) and
(\ref{eq:ortho_proj_B_bound}), it is clear that we only need to
provide a sharp upper bound for $\frac{1}{n}\sum_{i=1}^n\parallel
\Phi_i^T (I_M - B_* B_*^T) B\parallel_F^2$. Let $\overline{\Delta}
:= (I_M - B_* B_*^T) B$. Then from (\ref{eq:exp_U_2}), for $n$ large
enough,
\begin{equation}\label{eq:Delta_bar_bound}
\parallel \overline{\Delta} \parallel_F^2 \leq c_* \alpha_n^2
\end{equation}
for some $c_* > 0$. Then, using (\ref{eq:g_cond}), for all $i$,
\begin{eqnarray}\label{eq:ortho_proj_B_mean}
\mathbb{E}\parallel \Phi_i^T (I_M - B_* B_*^T) B \parallel_F^2 &=&
\sum_{j=1}^{m_i} \tr\left( \mathbb{E}
(\boldsymbol{\phi}(T_{ij})(\boldsymbol{\phi}(T_{ij}))^T)
\overline{\Delta}\overline{\Delta}^T\right)\nonumber\\
&\leq& c_{g,1} m_i \tr[\overline{\Delta}\overline{\Delta}^T] \leq
c_{g,1} \overline{m} \parallel \overline{\Delta} \parallel_F^2.
\end{eqnarray}
Combining (\ref{eq:ortho_proj_B_mean}) with
(\ref{eq:Delta_bar_bound}), (\ref{eq:Sigma_diff_upper1}) and
(\ref{eq:ortho_proj_B_bound}), we get, for sufficiently large $n$,
and some constant $C > 0$,
\begin{eqnarray*}\label{eq:Sigma_diff_upper3}
\frac{1}{n}\sum_{i=1}^n \mathbb{E} \parallel \Sigma_i -
\Sigma_{*i}\parallel_F^2 &\leq& C~ \overline{m}^2 \alpha_n^2
\end{eqnarray*}
Next, using (\ref{eq:ortho_proj_B_bound}), we have
\begin{eqnarray}\label{eq:ortho_proj_B_var}
\max_{1\leq i \leq n} \mbox{Var} (\parallel \Phi_i^T
\overline{\Delta} \parallel_F^2) &\leq& \max_{1\leq i \leq n}
\mathbb{E}
\parallel \Phi_i^T \overline{\Delta} \parallel_F^4 \nonumber\\
&\leq& c_{\phi,2} \overline{m} M \parallel \overline{\Delta}
\parallel_F^2
\mathbb{E} \parallel \Phi_i^T \overline{\Delta}
\parallel_F^2\nonumber\\
&\leq& c_{\phi,2} c_{g,1} \overline{m}^2 M \parallel \overline{\Delta}
\parallel_F^4 \nonumber\\
&\leq& C' M \overline{m}^2  \alpha_n^4 (1+\varepsilon_{n}) ~=:~
V_{2,n},
\end{eqnarray}
for some positive sequence $\varepsilon_{n} = o(1)$ and some
constant $C' > 0$, where in the last step we used
(\ref{eq:Delta_bar_bound}). Again, using Bernstein's inequality for
$\frac{1}{n}\sum_{i=1}^n\parallel \Phi_i^T \overline{\Delta}
\parallel_F^2$, which is a sum of independent variables bounded by
$K_{2,n} = c_{g,1} M \overline{m} \alpha_n^2(1+o(1))$, the result
follows (checking that, by (\ref{eq:M_condition_sparse}) or
(\ref{eq:m_M_condition}), we have, $K_{2,n} \sqrt{\frac{M\log n}{n}}
= o(\sqrt{V_{2,n}})$ and $\sqrt{V_{2,n}} \sqrt{\frac{M\log n}{n}} =
o(\overline{m}^2\alpha_n^2)$).

\vskip.15in\noindent{\bf Proof of Proposition 3 :} Write, $R_i =
\overline{\Sigma}_i^{1/2} (\Sigma_i^{-1} - \Sigma_{*i}^{-1})
\overline{\Sigma}_i^{1/2}$. We can bound $\parallel R_i\parallel_F$
as
\begin{eqnarray*}
\parallel R_i \parallel_F &\leq& \parallel
\overline{\Sigma}_i^{1/2} \Sigma_i^{-1/2}\parallel \parallel
\overline{\Sigma}_i^{1/2} \Sigma_{*i}^{-1/2}\parallel
\parallel \Sigma_i^{-1/2} \Sigma_{*i}^{1/2}\parallel \parallel
\Sigma_{*i}^{-1/2}(\Sigma_i -
\Sigma_{*i})\Sigma_{*i}^{-1/2}\parallel_F\nonumber\\
&\leq& \parallel \overline{\Sigma}_i^{1/2}
\Sigma_{*i}^{-1/2}\parallel^2
\parallel \Sigma_i^{-1/2} \Sigma_{*i}^{1/2}\parallel^2  \parallel
\Sigma_i - \Sigma_{*i}\parallel_F\nonumber\\
&\leq& (1+\parallel \overline{\Sigma}_i - \Sigma_{*i}
\parallel)
(1-\parallel \Sigma_i - \Sigma_{*i}\parallel)^{-1} \parallel
\Sigma_i - \Sigma_{*i}\parallel_F,
\end{eqnarray*}
where the third inequality is due to (\ref{eq:Sigma_eigen_bound}).
Note that, by condition {\bf C}, it follows that $\max_{1\leq i \leq
n}
\parallel \overline{\Sigma}_i - \Sigma_{*i}\parallel_F \leq C~
\overline{m}\overline{\beta}_n = o(1)$ for some constant $C
> 0$. Therefore, applying (\ref{eq:Sigma_diff_upper2}), we observe
that for $\mathbf{T} \in A_{\eta}$ with $A_{\eta}$ as in
(\ref{eq:T_eta_set}), for large enough $n$, $\parallel R_i
\parallel_F \leq 2 \parallel \Sigma_i - \Sigma_{*i}\parallel_F$.
Due to the Gaussianity of the observations, for any symmetric $m_i
\times m_i$ matrix $A$, the random variable  $\tr(A (S_i -
\overline{\Sigma}_i))$ has the same distribution as $\tr(D_i
(X_iX_i^T - I_{m_i}))$, where $D_i$ is a diagonal matrix of the
eigenvalues of $\overline{\Sigma}_i^{-1/2} A
\overline{\Sigma}_i^{-1/2}$, and $X_i \sim N(0,I_{m_i})$ are
independent. Therefore, using an exponential inequality for a
weighted sum of independent $\chi_1^2$ random variables, we have,
for ${\bf T} \in A_{\eta}$ and each $(B,\Lambda) \in
\Theta(\alpha_n)$,
\begin{eqnarray*}
&&\mathbb{P}\left(\left|\frac{1}{n}\sum_{i=1}^n
\tr((\Sigma_i^{-1}-\Sigma_{*i}^{-1})(S_i
-\overline{\Sigma}_i))\right| \leq d_{3,\eta} \sqrt{\frac{M \log
n}{n}} \left(\frac{1}{n} \sum_{i=1}^n
\parallel \Sigma_i - \Sigma_{*i}
\parallel_F^2\right)^{1/2} ~\Bigl|~ \mathbf{T}\right)\nonumber\\ &\geq&
1-O(n^{-(2+2\kappa)Mr -\eta}),
\end{eqnarray*}
for a constant $d_{3,\eta} >0$. Therefore, using
(\ref{eq:Sigma_diff_upper_concen}) we conclude the proof.

\vskip.15in\noindent{\bf Proof of Proposition 4 :} Using
Cauchy-Schwarz inequality twice, we can bound the last term in
(\ref{eq:loss_decomp_sparse}), which corresponds to \textit{model
bias}, as
\begin{eqnarray*}
\left|\frac{1}{2n}\sum_{i=1}^n \tr((\Sigma_i^{-1} -
\Sigma_{*i}^{-1})(\overline{\Sigma}_i - \Sigma_{*i}))\right| &\leq&
\frac{1}{2} \left[\frac{1}{n} \sum_{i=1}^n\parallel
\overline{\Sigma}_i -
\Sigma_{*i}\parallel_F^2\right]^{1/2}\left[\frac{1}{n} \sum_{i=1}^n
\parallel \Sigma_i^{-1} -
\Sigma_{*i}^{-1}\parallel_F^2\right]^{1/2}\nonumber\\
&\leq& \frac{1}{2} \max_{1\leq j \leq n} \parallel
\overline{\Sigma}_j - \Sigma_{*j}\parallel_F \left[\frac{1}{n}
\sum_{i=1}^n
\parallel \Sigma_i^{-1} -
\Sigma_{*i}^{-1}\parallel_F^2\right]^{1/2} \nonumber\\
&\leq& \frac{1}{2} \overline{m}\overline{\beta}_n \max_{1\leq j \leq
n} \parallel \Sigma_j^{-1} \parallel \parallel \Sigma_{*j}^{-1}
\parallel \left[\frac{1}{n}\sum_{i=1}^n \parallel \Sigma_i -
\Sigma_{*i}\parallel_F^2\right]^{1/2},
\end{eqnarray*}
where in the last step we used (\ref{eq:inverse_identity}). Thus,
the proof is finished by using condition {\bf C},
(\ref{eq:Sigma_eigen_bound}) and (\ref{eq:Sigma_diff_upper_concen}).

\vskip.1in\noindent Now, to the complete the proof of Theorem 1, we
give the details of the covering argument.

\vskip.15in\noindent{\bf Proof of Lemma 2 :} Using an expansion
analogous to (\ref{eq:loss_decomp_sparse}) and the upper bound in
(\ref{eq:Kullback_inequality_sparse}), and applying Cauchy-Schwarz
inequality, we have, for some constants $C_1, C_2  > 0$, on
$A_{4,\eta}$ and for $\mathbf{T} \in A_\eta$, for $n$ large enough,
\begin{eqnarray*}\label{eq:loss_bound_delta}
&&|L_n(B_1,\Lambda_1) -
L_n(B_2,\Lambda_2)|\nonumber\\
&\leq& \frac{1}{n}\sum_{i=1}^n\left[\parallel \Sigma_{1,i}^{-1} -
\Sigma_{2,i}^{-1} \parallel_F
\parallel S_i -\overline{\Sigma}_i\parallel_F +
C_1 \parallel \Sigma_{1,i} - \Sigma_{2,i} \parallel_F^2 +
\parallel \Sigma_{1,i}^{-1} - \Sigma_{2,i}^{-1} \parallel_F
\parallel \overline{\Sigma}_i - \Sigma_{2,i} \parallel_F \right]\nonumber\\
&\leq& \max_{1\leq i \leq n} \parallel \Sigma_{1,i}^{-1} -
\Sigma_{2,i}^{-1} \parallel_F \parallel \overline{\Sigma}_i\parallel
\max_{1\leq i \leq n} \parallel \overline{\Sigma}_i^{-1/2} S_i
\overline{\Sigma_i}^{-1/2} - I_{m_i}
\parallel_F + C_2 \max_{1\leq i \leq n}
\parallel \Sigma_{1,i} - \Sigma_{2,i} \parallel_F^2\nonumber\\
&& + \max_{1\leq i \leq n} \parallel \Sigma_{1,i}^{-1} -
\Sigma_{2,i}^{-1} \parallel_F (\max_{1\leq i \leq n} \parallel
\overline{\Sigma}_i - \Sigma_{*i}\parallel_F + \max_{1\leq i \leq n}
\parallel \Sigma_{2,i} - \Sigma_{*i}\parallel_F)
\nonumber\\
&\leq& d_{4,\eta} [\overline{m}^2\delta_n \overline{m} \log n +
\overline{m}^2 \delta_n^2 +
\overline{m}\delta_n(\overline{m}\overline{\beta}_n+
\sqrt{M}\overline{m} \alpha_n)] = o(\alpha_n^2).
\end{eqnarray*}
In the last step, we have used Lemma 4 (for the last term), the
identity (\ref{eq:inverse_identity}) in Appendix F, and the fact
that $\parallel \Sigma_{k,i}^{-1}\parallel \leq 1$ ($k=1,2$).

\vskip.15in\noindent{\bf Proof of Corollary 1:} The best rate
follows by direct calculation.

The near-optimality of the estimator requires proving that, for an
appropriately chosen subclass ${\cal C}$ of covariance kernels of
rank $r$, we have the following analog of Theorem 2 of Hall
\textit{et al.} (2006): for any estimator $\{\widehat
\psi_k\}_{k=1}^r$ of the eigenfunctions $\{\psi_k\}_{k=1}^r$, for
$n$ sufficiently large,
\begin{equation}\label{eq:rate_lower_bound}
\min_{1\leq k \leq r} \sup_{\overline{\Sigma}_0 \in {\cal C}}
\mathbb{E}\parallel \widehat \psi_k - \overline{\psi}_k\parallel_2^2
\geq C n^{-8/9},
\end{equation}
for some $C > 0$. Here the parameter space ${\cal C}$ consists of
covariance kernels of rank $r$ with eigenfunctions satisfying {\bf
A1}-{\bf A2}. Moreover, the random design satisfies {\bf B1}-{\bf
B2}, with $\overline{m}$ bounded above.

The derivation of the lower bound on the risk involves construction
of a finite, ``least favorable'' parameter set in ${\cal C}$ by
combining the constructions in Paul and Johnstone (2007) (for
obtaining lower bounds on risk in high-dimensional PCA) and Hall
\textit{et al.} (2006) (for functional data case). This construction
is as follows. Let $\phi_1^0,\ldots,\phi_r^0$ be a set of
orthonormal functions on $[0,1]$ which are four times continuously
differentiable, with fourth derivative bounded. Let $M_* \asymp
n^{1/9}$ be an integer appropriately chosen. Let
$\gamma_1,\ldots,\gamma_{M_*}$ be a set of basis functions that are
(i) orthonormal on $[0,1]$, and orthogonal to the set
$\{\phi_1^0,\ldots,\phi_r^0\}$; (ii) are four times continuously
differentiable and $\gamma_j$ is supported on an interval of length
$O(M_*^{-1})$ around the point $\frac{j}{M_*}$. One particular
choice for these functions is to let $\{\phi_k^0\}$ be the
translated periodized scaling functions of a wavelet basis at a
certain \textit{scale} with adequate degree of smoothness, and to
let $\{\gamma_j\}_{j=1}^{M_*}$ be the set of compactly supported,
orthonormal, periodized wavelet functions corresponding to the
scaling functions. Indeed, then we can choose $M_*$ to be an integer
power of 2. Note that, such a basis ($\{\phi_k^0:1\leq k \leq r\}
\cup \{\gamma_l : 1\leq l \leq M_*\}$) has the stability and
smoothness property commensurate with the orthonormalized B-spline
basis we are using for deriving the REML estimators. Next, let
$\overline{\lambda}_1
> \cdots > \overline{\lambda}_r > 0$ be fixed numbers satisfying
{\bf A1}. Finally, let us define a covariance kernel
$\overline{\Sigma}_0^{(0)}$ as
\begin{equation}\label{eq:Sigma_0_0}
\overline{\Sigma}_0^{(0)}(s,t) = \sum_{k=1}^r \overline{\lambda}_k
\phi_k^0(s) \phi_k^0(t),~~~s,t \in [0,1].
\end{equation}
Also, for each fixed $j$ in some index set ${\cal F}_0$ (to be
specified below), define
$$
[\psi_1^{(j)}(s):\cdots:\psi_r^{(j)}(s)] = \widetilde\Psi(s)
\overline{B}^{(j)}, ~~~s \in [0,1]
$$
where $\widetilde\Psi(s) =
(\phi_1^0(s),\ldots,\phi_r^0(s),\gamma_1(s),\ldots,\gamma_{M^*}(s))$
and $\overline{B}^{(j)}$ is an $(M_*+r)\times r$ matrix with
orthonormal columns (to be specified below). Then define
\begin{equation}\label{eq:Sigma_0_j}
\overline{\Sigma}_0^{(j)}(s,t) = \sum_{k=1}^r \overline{\lambda}_k
\psi_k^{(j)}(s)\psi_k^{(j)}(t),~~~~s,t \in [0,1],
\end{equation}
for $j \in {\cal F}_0$. We require that $\log |{\cal F}_0| \asymp
M_* \asymp n^{1/9}$, and $\parallel \overline{B}^{(j)} -
\overline{B}^{(j')}\parallel_F^2 \asymp n^{-8/9}$, for $j\neq j'$
and $j,j' \in {\cal F}_0 \cup\{0\}$. Here $\overline{B}^{(0)}$ is
the $(M_*+r)\times r$ matrix of basis coefficients of
$\overline{\Sigma}_0^{(0)}$ with columns $\mathbf{e}_k$, the $k$-th
canonical basis vector in $\mathbb{R}^{M_*+r}$.

The proof of the minimax lower bound is based on an application of
\textit{Fano's Lemma} (Yang and Barron, 1999), which requires
computation of the Kullback-Leibler divergence between two specific
values of the parameters. In order to apply \textit{Fano's lemma},
we need to choose ${\cal F}_0$ and $\overline{B}^{(j)}, j \in {\cal
F}_0$, such that
\begin{equation}\label{eq:Fano_ratio}
\frac{\mbox{ave}_{j \in {\cal F}_0} [\sum_{i=1}^n \mathbb{E}
K(\overline{\Sigma}_i^{(j)}, \overline{\Sigma}_i^{(0)})] + \log
2}{\log |{\cal F}_0|} \approx c \in (0,1),
\end{equation}
where $\overline{\Sigma}_i^{(j)}$ denotes the covariance of the
observation $i$ given $\{T_{il}\}_{l=1}^{m_i}$ under the model
parameterized by $\overline{\Sigma}_0^{(j)}$, and $\mathbb{E}$
denotes expectation with respect to the design points $\mathbf{T}$.
Under the assumptions on the design points, using the properties of
the basis functions $\{\phi_k^0\}_{=1}^r$ and
$\{\gamma_k\}_{k=1}^{M_*}$, and the computations carried out in the
proof of Proposition 2 (in Appendix B), in particular a nonrandom
bound analogous to the second bound appearing in Lemma 4, it is easy
to see that for $n$ large enough (so that $\parallel
\overline{B}^{(j)} - \overline{B}^{(0)}\parallel_F$ is sufficiently
small), we have
\begin{equation*}
\frac{1}{n} \sum_{i=1}^n K(\overline{\Sigma}_i^{(j)},
\overline{\Sigma}_i^{(0)}) \asymp \frac{1}{n} \sum_{i=1}^n
\parallel \overline{\Sigma}_i^{(j)} -
\overline{\Sigma}_i^{(0)}\parallel_F^2.
\end{equation*}
From this, and the property of the basis used to represent the
eigenfunctions, it follows that
\begin{equation}\label{eq:Fano_KL_bound}
\frac{1}{n} \sum_{i=1}^n \mathbb{E} K(\overline{\Sigma}_i^{(j)},
\overline{\Sigma}_i^{(0)}) \asymp \parallel \overline{B}^{(j)} -
\overline{B}^{(0)}\parallel_F^2.
\end{equation}
The task remains to construct ${\cal F}_0$ and $\overline{B}^{(j)}$
appropriately so that ${\cal C}_0 := \{\overline{\Sigma}_0^{(j)}:j
\in {\cal F}_0\}$ is in ${\cal C}$, for $n$ sufficiently large.

Following the proof of Theorem 2 in Paul and Johnstone (2007), we
first define $M_0 = [\frac{2 M_*}{9r}]$. Then define the $k$-th
column of $\overline{B}^{(j)}$ as
\begin{equation}\label{eq:B_bar_Fano}
\overline{B}_k^{(j)} = \sqrt{1 - \delta_k^2} \mathbf{e}_k + \delta_k
\sum_{l=1}^{M_*} z_{kl}^{(j)} \mathbf{e}_{r+l},~~~ k=1,\ldots,r,
\end{equation}
where $z_{kl}^{(j)}$ are appropriately chosen using a ``sphere
packing'' argument (to ensure that $\log |{\cal F}_0| \asymp M_*$),
and take values in $\{-M_0^{-1/2}, 0, M_0^{-1/2}\}$. Moreover, let
$S_k$ be the set of coordinates $l$ such that $z_{kl}^{(j)} \neq 0$
for some $j \in {\cal F}_0$. By construction, $S_k$ are disjoint for
different $k=1,\ldots,r$, and $|S_k| \sim M_*/r$. Hence,
\begin{equation}\label{eq:psi_k_j_repr}
\psi_k^{(j)} = \sqrt{1-\delta_k^2} \phi_k^0 + \delta_k \sum_{l \in
S_k} z_{kl}^{(j)} \gamma_l, ~~~k=1,\ldots,r.
\end{equation}
Furthermore, by the construction of $\{z_{lk}^{(j)}\}$, $\sum_{l\in
S_k} |z_{kl}^{(j)}|^2 = 1$, and for any $j\neq j'$ the vectors
$\mathbf{z}_k^{(j)} = (z_{kl}^{(j)})_{l \in S_k}$ and
$\mathbf{z}_k^{(j')} = (z_{kl}^{(j')})_{l \in S_k}$ satisfy
$\parallel \mathbf{z}_k^{(j)} - \mathbf{z}_k^{(j')}
\parallel_2\geq 1$. Therefore, from (\ref{eq:B_bar_Fano}) it follows
that the RHS of (\ref{eq:Fano_KL_bound}) is of the order
$\delta_k^2$, and hence in order that (\ref{eq:Fano_ratio}) is
satisfied, we need to choose $\delta_k \sim n^{-4/9} \asymp
M_*^{-4}$. It follows immediately from (\ref{eq:psi_k_j_repr}) that
(i) the eigenfunctions $\psi_1^{(j)},\ldots,\psi_r^{(j)}$ are
orthonormal, and four times continuously differentiable. Also, since
$\gamma_l$ is centered around $l/M_*$ with a support of the order
$O(M_*^{-1})$, it follows that, for only finitely many $l\neq l'$,
the support of $\gamma_{l'}$ overlaps with the support of
$\gamma_l$. Moreover, if $\gamma_l^{(s)}$ denotes the $s$-th
derivative of $\gamma_l$, then $\parallel \gamma_l^{(s)}
\parallel_\infty = O(M_*^{1/2 + s})$, for $s=0,1,\ldots,4$. Thus,
the choice $\delta_k \asymp M_*^{-4}$ ensures that, (ii) for each
$k=1,\ldots,r$, the fourth derivative of $\psi_k^{(j)}$  is bounded.
Hence, by appropriate choice of the constants, we have that ${\cal
C}_0 \subset {\cal C}$. Finally, arguing as in Paul and Johnstone
(2007), with an application of the \textit{Fano's Lemma} we conclude
(\ref{eq:rate_lower_bound}).

\section*{Appendix C : Proof of Proposition 5}

Using standard arguments, it can be shown that, given $\epsilon > 0$
sufficiently small (but fixed), we have constants $0 <
c_{1,\epsilon} < c_{2,\epsilon} < \infty$ such that for $\parallel
\Gamma_*^{-1/2}(\Gamma_* - \Gamma)\Gamma_*^{-1/2}
\parallel_F \leq \epsilon$,
\begin{equation}\label{eq:Kullback_inequality}
c_{1,\epsilon} \parallel \Gamma_*^{-1/2}(\Gamma_* -
\Gamma)\Gamma_*^{-1/2}
\parallel_F^2 \leq K(\Gamma,\Gamma_*) \leq c_{2,\epsilon}\parallel
\Gamma_*^{-1/2}(\Gamma_* - \Gamma)\Gamma_*^{-1/2} \parallel_F^2.
\end{equation}
Thus, it suffices to provide tight bounds for $\parallel
\Gamma_{*}^{-1/2}(\Gamma -\Gamma_{*})\Gamma_{*}^{-1/2}\parallel_F$.
We introduce some notations first. Define, $G = \frac{\sigma^2}{s}
\Lambda^{-1} + I_r$, $G_{*} = \frac{\sigma^2}{s} \Lambda_*^{-1} +
I_r$ and $\Delta = B \Lambda B^T - B_* \Lambda_* B_*^T$. Then,
\begin{equation*}
\Gamma^{-1} = \frac{1}{\sigma^2} (I_M - B (\frac{\sigma^2}{s}
\Lambda^{-1} + I_r)^{-1} B^T ) = \frac{1}{\sigma^2} (I_M - B G^{-1}
B^T),~~\mbox{and}~~\Gamma_{*}^{-1} = \frac{1}{\sigma^2} (I_M - B_*
G_*^{-1} B_*^T).
\end{equation*}
Moreover, due to {\bf A1'}, there exist constants, $c_3, c_4
> 0$, such that,
\begin{equation*}
c_3 \left(\frac{\sigma^2}{s}\right) \leq \sigma_{min}(I_r -
G_{*}^{-1}) \leq \sigma_{max}(I_r - G_{*}^{-1}) \leq c_4
\left(\frac{\sigma^2}{s}\right).
\end{equation*}
We express $I_M - B_* G_*^{-1} B_*^T$ as $(I_M - B_* B_*^T) + B_*
(I_r - G_*^{-1}) B_*^T$. Then we can express \\ $(\sigma^2/s)^2
\parallel \Gamma_{*}^{-1/2} (\Gamma - \Gamma_*)
\Gamma_{*}^{-1/2}\parallel_F^2 = \sigma^4
\parallel \Gamma_{*}^{-1/2}\Delta \Gamma_{*}^{-1/2}\parallel_F^2$ as
\begin{eqnarray}\label{eq:norm_bound_lower_1}
&& \tr[(I_M - B_* G_{*}^{-1} B_*^T) \Delta (I_M - B_* G_{*}^{-1}
B_*^T) \Delta ]\nonumber\\
&=& \tr[(I_M - B_* B_*^T) B \Lambda B^T (I_M - B_* B_*^T) B\Lambda
B^T] + 2\tr[B_* (I_r - G_{*}^{-1}) B_*^T B \Lambda B^T (I_M - B_*
B_*^T) B \Lambda B^T]\nonumber\\
&& + \tr[B_* (I_r - G_{*}^{-1}) B_*^T \Delta B_* (I_r - G_{*}^{-1})
B_*^T \Delta]\nonumber\\
&=& \parallel (I_M - B_* B_*^T) B\Lambda B^T (I_M - B_*
B_*^T)\parallel_F^2
\nonumber\\
&& + 2\tr[(I_r - G_{*}^{-1})^{1/2} B_*^T B \Lambda (B^T (I_M -
B_* B_*^T) B) \Lambda B^T B_* (I_r - G_{*}^{-1})^{1/2}]\nonumber\\
&& + \parallel (I_r - G_{*}^{-1})^{1/2} (B_*^T B \Lambda B^T B_* -
\Lambda_*) (I_r - G_{*}^{-1})^{1/2}\parallel_F^2\nonumber\\
&\geq& 2c_3 \lambda_r^2 (\sigma_{min}(B_*^T B))^2 \frac{\sigma^2}{s}
\parallel (I_M - B_* B_*^T) B\parallel_F^2 + c_3^2
\left(\frac{\sigma^2}{s}\right)^2 \parallel B_*^T B \Lambda B^T B_*
- \Lambda_*\parallel_F^2
\nonumber\\
&\geq& c_4 \left(\frac{\sigma^2}{s}\right)\parallel (I_M - B_*
B_*^T) B\parallel_F^2 + c_3^2 \left(\frac{\sigma^2}{s}\right)^2
\parallel B_*^T B \Lambda B^T B_* - \Lambda_*\parallel_F^2
\end{eqnarray}
for constants $c_3, c_4 > 0$. Now, since $(B,\Lambda) \in \widetilde
\Theta(\alpha_n)$, where $B = \mathbf{exp}(1,B_* A_U + C_U)$ it
follows that $\parallel A_U \parallel_F \leq \alpha_n
\sqrt{\frac{s}{\sigma^2}}$ and $\parallel C_U \parallel_F \leq
\alpha_n$. Moreover, from (\ref{eq:exp_U_1}), and using the fact
that $A_U = - A_U^T$, we have,
\begin{eqnarray*}
B_*^T B \Lambda B^T B_* - \Lambda_* &=& D + (A_U \Lambda - \Lambda
A_U) + O(\parallel A_U \parallel_F^2 + \parallel D \parallel_F^2 +
\parallel U \parallel_F (\parallel A_U \parallel_F + \parallel C_U
\parallel_F)).
\end{eqnarray*}
Since $A_U \Lambda - \Lambda A_U$ is symmetric, has zeros on the
diagonal, and its Frobenius norm is bounded below by $\min_{1\leq j
< k \leq r} (\lambda_j - \lambda_k) \parallel A_U \parallel_F$, and
$D$ is diagonal, it follows that for some constant $c_6 > 0$,
$$
\parallel B_*^T B \Lambda B^T B_* - \Lambda_* \parallel_F^2 \geq c_6
(\parallel D \parallel_F^2 + \parallel A_U \parallel_F^2) -
O((\frac{s}{\sigma^2})^{3/2}\alpha_n^3).
$$
From this, and using (\ref{eq:exp_U_2}) to approximate the first
term in (\ref{eq:norm_bound_lower_1}), it follows that for some
constant $c_7 > 0$,
\begin{eqnarray}\label{eq:Frob_norm_bound_lower_2}
&& \parallel (I_M - B_* G_{*i}^{-1} B_*^T)^{1/2} \Delta (I_M - B_*
G_{*i}^{-1} B_*^T)^{1/2}\parallel_F^2
\nonumber\\
&\geq& c_7 \left( \frac{\sigma^2}{s}\right) \left[\parallel
C_U\parallel_F^2 + \frac{\sigma^2}{s}  \parallel A_U\parallel_F^2 +
\frac{\sigma^2}{s}
\parallel D\parallel_F^2 - O(\sqrt{\frac{s}{\sigma^2}}\alpha_n^3)\right]\nonumber\\
&=& c_7 \left(\frac{\sigma^2}{s}\right) \alpha_n^2(1-o(1)).
\end{eqnarray}
The last equality is because $\alpha_n \sqrt{\frac{s}{\sigma^2}} =
o(1)$. Also, it is easy to show now that, for some $c_8 > 0$
\begin{equation}\label{eq:Frob_norm_bound_upper}
\parallel (I_M - B_* G_{*}^{-1} B_*^T)^{1/2} \Delta
(I_M - B_* G_{*}^{-1} B_*^T)^{1/2}\parallel_F^2
\leq c_8 \left(\frac{\sigma^2}{s}\right) \alpha_n^2(1+o(1)).
\end{equation}
Hence, from (\ref{eq:Frob_norm_bound_lower_2}) and
(\ref{eq:Frob_norm_bound_upper}), it follows that, there are
constants $c_9,c_{10} > 0$ such that, for sufficiently large $n$,
\begin{equation*}
c_9 \alpha_n \sqrt{\frac{s}{\sigma^2}} \leq \parallel
\Gamma_{*}^{-1/2}(\Gamma - \Gamma_{*}) \Gamma_{*}^{-1/2}\parallel_F
\leq c_{10} \alpha_n \sqrt{\frac{s}{\sigma^2}},
\end{equation*}
which, together with (\ref{eq:Kullback_inequality}), proves
(\ref{eq:Kullback_inequality_final}).

\section*{Appendix D : Score representation in the matrix case}

First, define the \textit{canonical metric} on the tangent space
${\cal T}_{B} \oplus \mathbb{R}^r$ of the parameter space
$\tilde{\Omega}=S_{M,r} \otimes \mathbb{R}^r$ (for $\theta =
(B,\zeta)$) by
$$
\langle X, Y \rangle_g = \langle X_B, Y_B \rangle_c + \langle
X_\zeta, Y_\zeta\rangle,~~\mbox{for}~~X_B,Y_B \in {\cal T}_B,
~X_{\zeta}, Y_{\zeta} \in \mathbb{R}^r,
$$
where $\langle X_B, Y_B \rangle_c = \tr(X_B^T (I_M - \frac{1}{2}
BB^T) Y_B)$ is the canonical metric on ${\cal S}_{M,r}$ and
$\langle X_\zeta, Y_\zeta\rangle = \tr(X_{\zeta}^T Y_{\zeta})$ is
the usual Euclidean metric. Next, for an arbitrary $\theta$, write
$\widetilde L_n(\theta) = F_n^1(\theta) + F_n^2(\theta)$, where
$F_n^1 (\theta) = \tr(\Gamma^{-1}\widetilde{S})$ and
$F_n^2(\theta) = \log |\Gamma| = \log|I_r + e^\zeta| = \log |I_r +
\Lambda|$. Similarly, we write $L(\theta;\theta_*) =
F^1(\theta;\theta_*) + F^2(\theta;\theta_*)$, where
$F^1(\theta;\theta_*) = \tr(\Gamma^{-1}\Gamma_*)$ and
$F^2(\theta;\theta_*) = F_n^2(\theta)$. Below, we shall only give
expressions for gradient and Hessian of $F_n^1(\cdot)$ and
$F_n^2(\cdot)$, since the gradient and Hessian of
$F^1(\cdot;\theta_*)$ and $F^2(\cdot;\theta_*)$ follow from these
(by replacing $\widetilde{S}$ with $\Gamma_*$).


\subsection*{Gradient and Hessian}

From Appendices B and D of Peng and Paul (2007), we obtain
expressions for the gradient and Hessian of $F_n^1(\cdot)$ and
$F_n^2(\cdot)$. We mainly follow the notations used there. Let, $P
:= P(\theta) = I_M + B \Lambda B^T$. Then $P^{-1} = I_M - B Q^{-1}
B^T$, where
\begin{equation*}
Q := Q(\theta) = \Lambda^{-1} + B^T B = \Lambda^{-1} + I_r
~~\Longrightarrow~~ Q^{-1} = \Lambda (I_r + \Lambda)^{-1}.
\end{equation*}
The fact that $Q$ is independent of $B$ is of importance in the
calculations throughout. Use $F_{n,B}^1(\cdot)$ to denote the
Euclidean gradient of $F_n^1(\cdot)$  w.r.t. $B$. It is easy to see
that $F_{n,B}^1(\theta) = -2\widetilde{S} B Q^{-1}$. Then, under the
canonical metric the intrinsic gradient is given by
\begin{equation*}
\nabla_B F_n^1(\theta) = F_{n,B}^1(\theta) - B (F_{n,B}^1(\theta))^T
B = 2[B Q^{-1} B^T \widetilde{S} B - \widetilde{S} B Q^{-1}].
\end{equation*}
Since $F_n^2(\theta)$ does not involve $B$, the Euclidean gradient
$F_{n,B}^2(\theta) = 0$, and hence $\nabla_B F_n^2(\theta) = 0$.
Therefore,
\begin{equation}\label{eq:bar_L_n_grad_B}
\nabla_B \widetilde L_n(\theta) = \nabla_B F_n^1(\theta) = 2[B
Q^{-1} B^T \widetilde{S} B - \widetilde{S} B Q^{-1}].
\end{equation}
Next, for $X_B \in {\cal T}_B$, let $G_{n,BB}^1(\cdot)(X_B)$ be
the Euclidean Hessian operator of $F_n^1(\cdot)$ evaluated at
$X_B$. It is computed as
\begin{equation*}
G_{n,BB}^1(\theta)(X_B) = -2 \widetilde{S} X_B Q^{-1}.
\end{equation*}
The Hessian operator of $\widetilde L_n(\cdot)$ w.r.t. $B$, equals
the Hessian operator of $F_n^1(\cdot)$ w.r.t. $B$. For $X_B, Y_B \in
{\cal T}_B$, it is given by
\begin{eqnarray}\label{eq:F_1_hessian}
H_{n,B}(\theta)(X_B,Y_B) &=& \tr\left(Y_B^T
G_{n,BB}^1(\theta)(X_B)\right) +
\frac{1}{2}\tr\left[\left((F_{n,B}^1(\theta))^T X_B B^T +
B^T X_B (F_{n,B}^1(\theta))^T\right)Y_B\right] \nonumber\\
&& - \frac{1}{2} \tr\left[\left(B^TF_{n,B}^1(\theta) +
(F_{n,B}^1(\theta))^T B\right)X_B^T (I_M-BB^T)Y_B\right].
\end{eqnarray}
For computing gradient and Hessian with respect to $\zeta$, we
only need to compute first and second derivatives of the function
$\widetilde L_n(\cdot)$ (equivalently, of $F_n^1(\cdot)$ and
$F_n^2(\cdot)$). Using calculations carried out in Appendix D of
Peng and Paul (2007), and the identity $P^{-1} B_k =
(1+\lambda_k)^{-1} B_k$ where $B_k$ is the $k$-th column of $B$,
$1\leq k \leq r$, we have,
\begin{equation}\label{eq:F_12_grad_zeta}
\frac{\partial F_n^1}{\partial \zeta_k}(\theta) = - e^{\zeta_k}
 B_k^T P^{-1}\widetilde{S}P^{-1}B_k =
 -\frac{\lambda_k}{(1+\lambda_k)^2} B_k^T \widetilde{S}
 B_k, ~~~\mbox{and}~~~
\frac{\partial F_n^2}{\partial \zeta_k}(\theta) = e^{\zeta_k} B_k^T
P^{-1} B_k = \frac{\lambda_k}{1+\lambda_k}.
\end{equation}
Thus
\begin{equation*}
\nabla_\zeta \widetilde L_n(\theta) =
\mbox{diag}\left(\frac{\lambda_k}{(1+\lambda_k)^2}(1+\lambda_k -
B_k^T \widetilde{S}B_k)\right)_{k=1}^r.
\end{equation*}
Since $B_k^T P^{-1} B_l = 0$, for $1\leq k \neq l\leq r$, it follows
that $\frac{\partial^2 F_n^i}{\partial \zeta_k
\partial \zeta_l}(\theta) = 0$ for $k \neq l$, $i=1,2$. Also,
\begin{equation*}\label{eq:F_1_Hessian_zeta}
\frac{\partial^2 F_n^1}{\partial\zeta_k^2}(\theta) = e^{\zeta_k}
B_k^T P^{-1} \widetilde{S}P^{-1} B_k \left[2e^{\zeta_k} (B_k^T
P^{-1} B_k) - 1\right] = \frac{\lambda_k(\lambda_k
-1)}{(1+\lambda_k)^3} B_k^T \widetilde{S} B_k,
\end{equation*}
\begin{equation*}\label{eq:F_2_Hessian_zeta}
\frac{\partial^2 F_n^2}{\partial\zeta_k^2}(\theta) = e^{\zeta_k}
B_k^T P^{-1} B_k \left[1-e^{\zeta_k} (B_k^T P^{-1} B_k)\right] =
\frac{\lambda_k}{(1+\lambda_k)^2}.
\end{equation*}
Thus, the Hessian operator of $\widetilde L_n(\cdot)$ w.r.t. $\zeta$
is given by
\begin{equation}\label{eq:L_bar_n_Hessian_zeta}
H_{n,\zeta}(\theta) =
\mbox{diag}\left(\frac{\lambda_k}{(1+\lambda_k)^3}
\left((\lambda_k -1)B_k^T \widetilde{S} B_k +
(1+\lambda_k)\right)\right)_{k=1}^r.
\end{equation}

\subsection*{Boundedness and inversion of
$H(\theta_*;\theta_*)$}

As discussed in Section \ref{sec:proof_thm4}, the Hessian operator
$H(\theta_*;\theta_*)$ is ``block diagonal''. So, we only need to
show the boundedness and calculate the inverse of
$H_B(\theta_*;\theta_*)$ and $H_\zeta(\theta_*;\theta_*)$. First
note that, from (\ref{eq:L_bar_n_Hessian_zeta}) we have,
$$
H_{\zeta}(\theta_*;\theta_*) =
\mbox{diag}\left(\frac{\lambda_{*k}}{(1+\lambda_{*k})^3}
\left((\lambda_{*k} -1)B_{*k}^T \Gamma_* B_{*k} +
(1+\lambda_{*k})\right)\right)_{k=1}^r  =
\Lambda_*^2(I_r+\Lambda_*)^{-2},
$$
which is clearly positive definite with eigenvalues bounded away
from $0$ and $\infty$, due to conditions {\bf A1'} and {\bf C''}.

Next, we show that $H_{B}(\theta_*;\theta_*)(X,X) \geq C \langle X,X
\rangle_c$, for some $C > 0$, for all $X \in {\cal T}_{B_*}$.
Define
$F_B^1(\theta_*;\theta_*) = \mathbb{E}_{\theta_*}
F_{n,B}^1(\theta_*)$ and $G_{BB}^1(\theta_*;\theta_*) =
\mathbb{E}_{\theta_*} G_{n,BB}^1(\theta_*)$. Note that
$H_{B}(\theta_*;\theta_*)$ is obtained by replacing
$F_{n,B}^1(\theta_*)$ and $G_{n,BB}^1(\theta_*)$  by
$F_B^1(\theta_*)$ and $G_{BB}^1(\theta_*;\theta_*)$, respectively,
in (\ref{eq:F_1_hessian}). Observe that,
$$
F_B^1(\theta_*) = -2\Gamma_* B_* Q_*^{-1} = -2 B_* (I_r + \Lambda_*)
Q_*^{-1} = - 2B_* \Lambda_*,
$$
where $Q_* = Q(\theta_*) = \Lambda_*^{-1}(I_r+\Lambda_*)$, and we
have used the fact that $\Gamma_* B_* = B_* (I_r + \Lambda_*)$. For
notational simplicity we use $F_B^1$ to denote $F_B^1(\theta_*)$.
Note that, for $X \in {\cal T}_{B_*}$, $X = B_* A_X + (I- B_*B_*^T)
C_X$, where $A_X = - A_X^T \in \mathbb{R}^{r\times r}$, and $C_X \in
\mathbb{R}^{M\times r}$. Using this representation, for any $X,Y \in
{\cal T}_{B_*}$, we have
\begin{eqnarray}\label{eq:trace_term_1}
&&\frac{1}{2}\tr\left[\left((F_{B}^1)^T X B_*^T + B_*^T X
(F_{B}^1)^T\right)Y\right] ~=~
-\tr\left[\Lambda_* B_*^T X B_*^T Y + B_*^T X \Lambda_* B_*^T Y\right]\nonumber\\
&=& \tr\left[\Lambda_* X^T B_* B_*^T Y + B_*^T X \Lambda_* Y^T
B_*\right] ~=~2\tr\left[\Lambda_* X^T B_* B_*^T Y\right],
\end{eqnarray}
and
\begin{eqnarray}\label{eq:trace_term_2}
&&-\frac{1}{2} \tr\left[\left(B_*^T(F_{B}^1) + (F_{B}^1)^T
B_*\right)X^T \left(I_M - B_*B_*^T\right)Y\right]\nonumber\\
&=& \tr\left[ (B_*^T B_* \Lambda_*  + \Lambda_* B_*^T B_*) X^T
(I_M-B_*B_*^T)Y\right] ~=~ 2\tr\left[\Lambda_* X^T (I_M-B_*
B_*^T)Y\right].
\end{eqnarray}
Next, notice that, for $X \in {\cal T}_{B_*}$,
$G_{BB}^1(\theta_*;\theta_*)(X) = -2 \Gamma_* X Q_*^{-1}$.
Therefore,
\begin{eqnarray}\label{eq:trace_term_3}
\tr\left[Y^T G_{BB}^1(\theta_*;\theta_*)(X)\right] &\hskip-.1in=&
\hskip-.1in
-2\tr\left[Y^T \Gamma_* X Q_*^{-1}\right] \nonumber\\
&\hskip-.1in=& \hskip-.1in -2\tr\left[Y^T (I_M - B_*B_*^T) X
Q_*^{-1}] -2 \tr[Y^T B_*(I_r+\Lambda)B_*^T X Q_*^{-1}\right].
\end{eqnarray}
Now, combining (\ref{eq:trace_term_1}), (\ref{eq:trace_term_2}) and
(\ref{eq:trace_term_3}), and using the definition of
$H_{B}(\theta_*;\theta_*)$, and the facts that $I_r - Q_*^{-1} =
(I_r + \Lambda_*)^{-1}$, $X^T B_* = A_X^T = - A_X$ and $B_*^T Y =
A_Y = - A_Y^T$, after some simple algebra we have,
\begin{eqnarray}\label{eq:Hessian_repr_1}
&&H_{B}(\theta_*;\theta_*)(X,Y)
\nonumber\\
&=& 2\left(\tr\left[X^T B_*\Lambda_* B_*^T Y
(I_r+\Lambda_*)^{-1}\right] - \tr\left[X^T B_* B_*^T Y
Q_*^{-1}\right]\right) + 2\tr\left[\Lambda_*^2(I_r+\Lambda_*)^{-1}
X^T (I_r - B_* B_*^T) Y\right]
\nonumber\\
&&
\end{eqnarray}
Again, since $A_X = - A_X^T$ and $A_Y = - A_Y^T$, denoting by
$A_{X,ij}$ and $A_{Y,ij}$, the $(i,j)$-th element of $A_X$ and $A_Y$
respectively, we have
\begin{eqnarray}\label{eq:diff_term_Hessian}
&& \tr\left[X^T B_*\Lambda_* B_*^T Y (I+\Lambda_*)^{-1}\right] -
\tr\left[X^T B_*
B_*^T Y Q_*^{-1}\right]\nonumber\\
&=& - \tr\left[A_X \left(\Lambda_* A_Y (I_r+\Lambda_*) - A_Y
\Lambda_*(I_r+\Lambda_*)^{-1}\right)\right]\nonumber\\
&=& -\sum_{i=1}^r\sum_{j=1}^r A_{X,ij}
\left(\frac{\lambda_{*j}}{1+\lambda_{*i}} A_{Y,ji} -
\frac{\lambda_{*i}}{1+\lambda_{*i}} A_{Y,ji}\right) \nonumber\\
&=& \sum_{i=1}^r\sum_{j=1}^r A_{X,ij}A_{Y,ij}
\left(\frac{\lambda_{*j} - \lambda_{*i}}{1+\lambda_{*i}}\right) ~=~
\sum_{i=1}^{j-1}\sum_{j=i+1}^r A_{X,ij}A_{Y,ij} \left[(\lambda_{*i}
- \lambda_{*j}) \left(\frac{1}{1+\lambda_{*j}} -
\frac{1}{1+\lambda_{*i}}\right)\right]\nonumber\\
&=& \sum_{i=1}^{j-1}\sum_{j=i+1}^r A_{X,ij}A_{Y,ij}
\frac{(\lambda_{*i} -
\lambda_{*j})^2}{(1+\lambda_{*i})(1+\lambda_{*j})} = \frac{1}{2}
\sum_{i=1}^r\sum_{j=1}^r A_{X,ij}A_{Y,ij} \frac{(\lambda_{*i} -
\lambda_{*j})^2}{(1+\lambda_{*i})(1+\lambda_{*j})}.
\end{eqnarray}
Since $\min_{1\leq k \neq k' \leq r} (\lambda_{*k} -
\lambda_{*k'})^2 (1+\lambda_{*k})^{-2} \geq C_{1*}$, and
$\lambda_{*r} \geq  C_{*2}$, for some constants $C_{1*},C_{*2} >0$
(value depending on $\overline{c}_1$ and $\overline{c}_2$ appearing
in {\bf A1'}), it follows from (\ref{eq:Hessian_repr_1}) and
(\ref{eq:diff_term_Hessian}) that for $X \in {\cal T}_{B_*}$,
\begin{eqnarray*}\label{eq:Hessian_lower_bound}
H_{B}(\theta_*;\theta_*)(X,X) &\geq& C_{*1} \tr\left(X^T B_* B_*^T
X\right) +
2C_{*2}\tr\left(X^T (I_M - B_* B_*^T) X\right) \nonumber\\
&\geq& C_{*3}\tr\left(X^T X\right),
\end{eqnarray*}
where $C_{*3} = \min\{C_{*1}, 2C_{*2}\}$. This proves that
$H_{B}(\theta_*;\theta_*)(X,X)$ is bounded below in the Euclidean
norm and hence in the canonical metric because of the norm
equivalence. An upper bound follows similarly.


\vskip.15in\noindent{\bf Proof of Corollary 3 :} From
(\ref{eq:Hessian_repr_1}) and (\ref{eq:diff_term_Hessian}), we can
derive an explicit expression of $H_B^{-1}(\theta_*;\theta_*)$. Note
that $H_B^{-1}(\theta_*;\theta_*)(X)$ is defined as
$$
H_B\left(\theta_*;\theta_*)(H_B^{-1}(\theta_*;\theta_*)(X),Y\right)
= \langle X,Y\rangle_c,~~~\mbox{for any}~Y \in {\cal T}_{B_*}.
$$
Therefore, for $X \in {\cal T}_{B*}$ and $A_X = B_*^T X$,
\begin{eqnarray*}\label{eq:Hessian_B_L_inverse}
H_B^{-1}(\theta_*;\theta_*)(X) &=&  \frac{1}{2} B_*
\left(\left(\frac{(1+\lambda_{*i})(1+\lambda_{*j})}{(\lambda_{*i}
- \lambda_{*j})^2} A_{X,ij}\right)\right) + \frac{1}{2} (I_M -
B_*B_*^T) X\Lambda_*^{-2}(I_r+\Lambda_*).
\end{eqnarray*}
Using this, we can now get an explicit expression for
$H_B^{-1}(\theta_*;\theta_*)(\nabla_B \widetilde L_n(\theta_*))$.
From (\ref{eq:bar_L_n_grad_B}), we have
\begin{eqnarray*}\label{eq:B_*_grad_B}
B_*^T \nabla_B \widetilde L_n(\theta_*) &=& 2 \left[Q_*^{-1} B_*^T
\widetilde{S} B_* - B_*^T \widetilde{S} B_* Q_*^{-1}\right] =
2\left(\left( B_{*i}^T \widetilde{S} B_{*j}
\Bigl(\frac{\lambda_{*i}}{1+\lambda_{*i}}-\frac{\lambda_{*j}}{1+\lambda_{*j}}\Bigr)\right)\right)
\nonumber\\
&=& 2\left(\left( \frac{(\lambda_{*i} -
\lambda_{*j})}{(1+\lambda_{*i})(1+\lambda_{*j})}B_{*i}^T
\widetilde{S} B_{*j}\right)\right).
\end{eqnarray*}
Also,
\begin{equation*}\label{eq:Pi_*_grad_B}
(I_M - B_* B_*^T) \nabla_B \widetilde L_n(B_*,\Lambda_*)  = - 2 (I_M
- B_* B_*^T) \widetilde{S} B_* Q_*^{-1}.
\end{equation*}
Thus, it follows that
\begin{eqnarray*}
 - H_B^{-1}(\theta_*;\theta_*)\left(\nabla_B
\widetilde L_n(\theta_*)\right) &=& -\left(\left(
\frac{1}{(\lambda_{*i} - \lambda_{*j})}B_{*i}^T \widetilde{S}
B_{*j}\right)\right) + (I_M - B_*
B_*^T) \widetilde{S} B_* \Lambda_*^{-1}\nonumber\\
&=& -\left[\mathbf{R}_1 \widetilde{S} B_{*1} : \cdots :
\mathbf{R}_r \widetilde{S} B_{*r}\right].
\end{eqnarray*}

\section*{Appendix E : Gradient and Hessian on product manifolds}

In this section, we give a brief outline of the intrinsic geometry
associated with the product manifold of two Riemannian manifolds,
and as an application we consider the manifold ${\cal S}_{M,r}
\otimes \mathbb{R}^r$, which is the parameter space for
$(B,\zeta)$ in our problem.

\subsection*{Product Manifolds}

Consider two Riemannian manifolds: $\mathcal{M},  \mathcal{N}$ with
metrics $g_M$ and $g_N$, respectively. The product manifold
$\mathcal{P}$ of $\mathcal{M}, \mathcal{N}$ is then defined as:
$$
\mathcal{P} :=\mathcal{M} \otimes \mathcal{N}=\left\{(x,y): x \in
\mathcal{M}, y \in \mathcal{N}\right\}
$$
with the tangent space at a point $p=(x,y) \in \mathcal{P}$,
$$
\mathcal{T}_p\mathcal{P} :=\mathcal{T}_x\mathcal{M}\oplus
\mathcal{T}_y\mathcal{N}
$$
where $\mathcal{T}_x\mathcal{M}, \mathcal{T}_y\mathcal{N}$ are
tangent spaces of $\mathcal{M}, \mathcal{N}$ at points $x,y$,
respectively. The Riemannian metric $g$ on the tangent space
$\mathcal{T}\mathcal{P}$ is naturally defined as
$$
\langle T_1, T_2 \rangle_g := \langle \xi_1, \xi_2
\rangle_{g_M}+\langle \eta_1, \eta_2 \rangle_{g_N},
$$
where $T_i=(\xi_i,\eta_i) \in \mathcal{T}\mathcal{P}$, with $\xi_i
\in \mathcal{T}\mathcal{M}$ and $\eta_i \in \mathcal{T}\mathcal{N}$
($i=1,2$).

By the above definition of the product manifold $\mathcal{P}$, the
intrinsic gradient and Hessian of a smooth function $f$ defined on
$\mathcal{P}$ are as follows:
\begin{itemize}
\item Gradient:
\begin{equation*}\label{eqn:prod_grad}
 \nabla f = (\nabla_\mathcal{M}f_{\mathcal{M}},
\nabla_\mathcal{N} f_{\mathcal{N}}),
\end{equation*}
where $f_{\mathcal{M}}$ ($f_{\mathcal{N}}$) is $f$ viewed as a
function on $\mathcal{M}$ ($\mathcal{N}$); and
$\nabla_{\mathcal{M}}$ ($\nabla_{\mathcal{N}}$ ) denotes the
gradient operator for functions defined on $\mathcal{M}$
($\mathcal{N}$).

\item
Hessian: for $T_i=(\xi_i, \eta_i) \in \mathcal{T}\mathcal{P}$
($i=1,2$),
\begin{eqnarray*} \label{eqn:prod_hess}
H_f(T_1,T_2)&=&H_{f_\mathcal{M}}(\xi_1,\xi_2)+\langle
\nabla_{\mathcal{N}} \langle \nabla_{\mathcal{M}} f_{\mathcal{M}},
\xi_1 \rangle_{g_M}, \eta_2 \rangle_{g_N} \nonumber\\
&+&\langle \nabla_{\mathcal{M}} \langle \nabla_{\mathcal{N}}
f_{\mathcal{N}}, \eta_1 \rangle_{g_N}, \xi_2
\rangle_{g_M}+H_{f_\mathcal{N}}(\eta_1,\eta_2).
\end{eqnarray*}
The above expression is derived from the bi-linearity of the Hessian
operator and its definition. Also note that
$$
\langle \nabla_{\mathcal{N}} \langle \nabla_{\mathcal{M}}
f_{\mathcal{M}}, \xi_1 \rangle_{g_M}, \eta_2 \rangle_{g_N}=\langle
\nabla_{\mathcal{M}} \langle \nabla_{\mathcal{N}} f_{\mathcal{N}},
\eta_2 \rangle_{g_N}, \xi_1 \rangle_{g_M}.
$$
\end{itemize}

\subsection*{Application to the product of a Stiefel manifold
and an Euclidean space}

Consider the special case: $\mathcal{M}=\mathcal{S}_{M,r}$ with the
canonical metric $\langle \cdot,\cdot \rangle _c$, and
$\mathcal{N}=\mathbb{R}^d$ with Euclidean metric. For a point
$p=(B,x)$ on the product manifold $\mathcal{P}$, the tangent space
is
$$
\mathcal{T}_p\mathcal{P}=\mathcal{T}_B \mathcal{M} \oplus
\mathcal{T}_x \mathcal{N},
$$
where
$$
{\cal T}_B {\cal M} = \{ \mathit{\Delta} \in \mathbb{R}^{M\times r}
: B^T \mathit{\Delta} = - \mathit{\Delta}^T B \},
\qquad\mbox{and}\qquad  \mathcal{T}_x \mathcal{N}=\mathbb{R}^d.
$$
For a smooth function $f$ defined on the product space
$\mathcal{P}$:
\begin{itemize}
\item Gradient (at $p$):
\begin{eqnarray*}
\nabla f |_p=\left(\nabla_{\mathcal{M}} f, \frac{\partial
f}{\partial x}\right) \Bigl|_p,
\end{eqnarray*}
where $\nabla_{\mathcal{M}} f\left|_{p}\right. = f_B - B f_B^T B$
(with $f_B=\frac{\partial f}{\partial B}$).

\item
Hessian operator (at $p$): for $T=(\Delta,a)$, and for$X =(X_B,\eta)
\in \mathcal{T}_p\mathcal{P}$,
\begin{eqnarray}\label{eq:prod_Hessian_operator}
H_f(T,X)|_p=H_{f_{\mathcal{M}}}(\Delta,X_B)+ \langle \frac{\partial
}{\partial x}\langle \nabla_{\mathcal{M}}f, \Delta
 \rangle_c,\eta \rangle+\langle
\frac{\partial }{\partial x}\langle \nabla_{\mathcal{M}}f, X_B
 \rangle_c,a \rangle+a^T \frac{\partial^2 f}{\partial x^2}\eta,
\end{eqnarray}
where
$$
H_{f_{\mathcal{M}}}(\mathit{\Delta},X_B)\left|_{p}\right. =
f_{BB}(\mathit{\Delta},X_B) + \frac{1}{2} Tr\left[(f_B^T
\mathit{\Delta} B^T + B^T \mathit{\Delta} f_B^T)X_B\right] -
\frac{1}{2} Tr\left[ (B^T f_B + f_B^T B)
\mathit{\Delta}^T\mathit{\Pi} X_B\right],$$ with $\mathit{\Pi} = I -
BB^T$.

\item Inverse of Hessian operator (at $p$):  for $G \in
\mathcal{T}_p\mathcal{P}$, $T=H_f^{-1}(G)|_p$ is defined as:
$T=(\Delta, a) \in \mathcal{T}_p\mathcal{P}$ such that for any
$X=(X_B, \eta) \in \mathcal{T}_p\mathcal{P}$ the following
equation is satisfied
\begin{equation*}\label{eq:Hessian_prod_star}
H_f(T,X)|_p=\langle G, X \rangle_g \quad .
\end{equation*}

\end{itemize}

\section*{Appendix F : Some inequalities involving matrices}

In this paper we make frequent use of the following matrix
inequalities:
\begin{itemize}
\item For any $A$, $B$,
\begin{equation*}\label{eq:norm_inequality}
\parallel AB \parallel_F \leq \parallel A
\parallel_F
\parallel B
\parallel, ~~~~\mbox{and}~~~\parallel AB \parallel_F \geq \parallel A \parallel_F
\lambda_{min}(B), ~~~(\mbox{for}~~B~~ \mbox{positive definite})
\end{equation*}
where $\lambda_{min}(B)$ is the smallest eigenvalue of $B$. Also, if
$A$ and $B$ are invertible then
\begin{equation}\label{eq:inverse_identity}
A^{-1} - B^{-1} = A^{-1}(B - A)B^{-1} = B^{-1}(B - A)A^{-1}.
\end{equation}

\item \textit{Weilandt's inequality (Horn and Johnson (1994)):}
For symmetric $p \times p$
matrices $A$, $B$ with eigenvalue sequences $\lambda_1(A) \geq
\cdots \geq \lambda_p(A)$ and $\lambda_1(B) \geq \cdots \geq
\lambda_p(B)$, respectively,
\begin{equation}\label{eq:weilandt}
\sum_{i=1}^p|\lambda_i(A) - \lambda_i(B)|^2 \leq \parallel A - B
\parallel_F^2
\end{equation}
\item \textit{Eigenvector perturbation (Paul (2005)):}
Let $A$ be a $p\times p$ positive semidefinite matrix, with $j$-th
largest eigenvalue $\lambda_j(A)$ with corresponding eigenvector
$\mathbf{p}_j$, and $\tau_j := \max\{(\lambda_{j-1}(A) -
\lambda_j(A))^{-1},(\lambda_j(A)- \lambda_{j+1}(A))^{-1}\}$ is
bounded (we take $\lambda_0(A) = \infty$ and $\lambda_{p+1}(A) =
0$). Let $B$ be a symmetric matrix. If $\mathbf{q}_j$ denotes the
eigenvector of $A+B$ corresponding to the $j$-th largest eigenvalue
(which is of multiplicity 1, for $\parallel B \parallel$ small
enough, by (\ref{eq:weilandt})), then (assuming without loss of
generality $\mathbf{q}_j^T \mathbf{p}_j > 0$),
\begin{equation}\label{eq:eigenvec_perturb}
\parallel \mathbf{q}_j - \mathbf{p}_j\parallel \leq 5~
\frac{\parallel B \parallel}{\tau_j} + 4 \left(\frac{\parallel B
\parallel}{\tau_j}\right)^2.
\end{equation}
\end{itemize}


\section*{Reference}

\begin{enumerate}
\item
Antoniadis, A. and Sapatinas, T. (2007) : Estimation and inference
in functional mixed-effects models. \textit{Computational Statistics
and Data Analysis} {\bf 51}, 4793-4813.

\item
Bickel, P. J.  and Levina, E. (2007) : Covariance regularization by
thresholding.  \textit{Technical report  \#744},  Department of
Statistics, UC Berkeley. \textit{Annals of Statistics}. To appear.

\item
Bickel, P. J.  and Levina, E. (2008) : Regularized estimation of
large covariance matrices.  \textit{Annals of Statistics} {\bf 36},
199-227.

\item
Burman, P. (1985) : A data dependent approach to density estimation.
\textit{Zeitschrift f\"{u}r Wahrscheinlichkeitstheorie und verwandte
Gebiete} {\bf 69}, 609-628.

\item
Chen, A. and Bickel, P. J. (2006) : Efficient independent component
analysis. \textit{Annals of Statistics} {\bf 34}, 2825-2855.

\item
de Boor, C. (1974) : Bounding the error in spline interpolation.
\textit{SIAM Review} {\bf 16}, 531-544.

\item
de Boor, C. (1978) : \textit{A Practical Guide to Splines}.
Springer-Verlag.

\item
deVore, R. A. and Lorentz, G. G. (1993) : \textit{Constructive
Approximation}. Springer.

\item
El Karoui, N. (2007) : Operator norm consistent estimation of large
dimensional sparse covariance matrices. \textit{Annals of
Statistics}. To appear.

\item
Davidson, K. R., and Szarek, S. (2001) : Local operator theory,
random matrices and Banach spaces. In \textit{``Handbook on the
Geometry of Banach spaces''}, ({\bf V. 1}, Johnson, W. B.,
Lendenstrauss, J. eds.), 317-366, Elsevier Science.

\item
Edelman, A., Arias, T. A. and Smith, S. T. (1998) : The geometry of
algorithms with orthogonality constraints, \textit{SIAM Journal on
Matrix Analysis and Applications} {\bf 20}, 303-353.

\item
Fan, J., Fan, Y. and Lv., J. (2006) : High dimensional covariance
matrix estimation using a factor model. \textit{Journal of
Econometrics}, (to appear).

\item
Ferraty, F. and Vieu, P. (2006) : \textit{Nonparametric Functional
Data Analysis : Theory and Practice}. Springer.

\item
Hall, P. and Hosseini-Nasab, M. (2006) : On properties of functional
principal components analysis. \textit{Journal of the Royal
Statistical Society, Series B} {\bf 68}, 109-126.

\item
Hall, P., M\"{u}ller, H.-G. and Wang, J.-L. (2006) Properties of
principal component methods for functional and longitudinal data
analysis. \textit{Annals of Statistics} {\bf 34}, 1493-1517.

\item
Horn, R. A. and Johnson, C. R. (1994) : \textit{Topics in Matrix
Analysis}. Cambridge University Press.

\item
James, G. M., Hastie, T. J. and Sugar, C. A. (2000) : Principal
component models for sparse functional data. \textit{Biometrika},
{\bf 87}, 587-602.


\item
Kato, T. (1980) : \textit{Perturbation Theory of Linear Operators}.
Springer-Verlag.

\item
Muirhead, R. J. (1982) : \textit{Aspects of Multivariate Statistical
Theory}, John Wiley \& Sons.

\item
Oller, J. M. and Corcuera, J. M. (1995) : Intrinsic analysis of
statistical estimation. \textit{Annals of Statistics} {\bf 23},
1562-1581.

\item
Paul, D. (2005) : \textit{Nonparametric Estimation of Principal
Components}. Ph. D. Thesis. Stanford University.

\item Paul, D. and Johnstone, I. M. (2007) : Augmented sparse
principal component analysis for high dimensional data.
\textit{Working Paper}.\\
(\texttt{http://anson.ucdavis.edu/$\sim$debashis/techrep/augmented-spca.pdf})

\item Peng, J. and Paul, D. (2007) : A
geometric approach to maximum likelihood estimation of covariance
kernel from sparse irregular longitudinal data. \textit{Technical
Report}. arXiv:0710.5343v1 [stat.ME]. (also at
\texttt{http://anson.ucdavis.edu/$\sim$jie/pd-cov-likelihood-technical.pdf})

\item
Ramsay, J. and Silverman, B. W. (2005) : \textit{Functional Data
Analysis, 2nd Edition}. Springer.


\item
Yang, A. and Barron, A. (1999) : Information-theoretic determination
of minimax rates of convergence. \textit{Annals of Statistics} {\bf
27}, 1564-1599.

\item
Yao, F., M\"{u}ller, H.-G. and Wang, J.-L. (2005) : Functional data
analysis for sparse longitudinal data. \textit{Journal of the
American Statistical Association} {\bf 100}, 577-590.

\end{enumerate}
D. Paul\\
Department of Statistics\\
University of California\\
Davis, CA 95616\\
debashis@wald.ucdavis.edu\\

\vskip 0.1in

\noindent J. Peng\\
Department of Statistics\\
University of California\\
Davis, CA 95616\\
jie@wald.ucdavis.edu\\

\end{document}